\documentclass{article}
\usepackage{latexsym}
\usepackage{amssymb}
\begin{document}
\def\P{\Bbb P}
\def\Q{\Bbb Q}
\title{Forcings constructed along morasses}
\author{Bernhard Irrgang}
\date{}
\maketitle

\begin{abstract}
We further develop a previously introduced method of constructing forcing notions with the help of morasses. There are two new results: (1) If there is a simplified $(\omega_1,1)$-morass, then there exists a ccc forcing of size $\omega_1$ that adds an $\omega_2$-Suslin tree.   (2) If there is a simplified $(\omega_1,2)$-morass, then there exists a ccc forcing of size $\omega_1$ that adds  a $0$-dimensional Hausdorff topology $\tau$ on $\omega_3$ which has spread $s(\tau)=\omega_1$. While (2) is the main result of the paper, (1) is only an improvement of a previous result, which is based on a simple observation. Both forcings preserve GCH. To show that the method can be changed to produce models where CH fails, we give an alternative construction of Koszmider's model in which there is a chain $\langle X_\alpha \mid \alpha < \omega_2\rangle$ such that $X_\alpha \subseteq \omega_1$, $X_\beta -X_\alpha$ is finite and $X_\alpha-X_\beta$ has size $\omega_1$ for all $\beta < \alpha < \omega_2$. 
\end{abstract}
\section{Introduction}
In a previous paper \cite{Irrgang4}, we introduced a method of constructing a forcing along a simplified $(\kappa,1)$-morass such that the forcing satisfies a chain condition. The basic idea is simple: We try to generalize iterated forcing with finite support (FS). Classical iterated forcing with finite support as introduced by Solovay and Tennenbaum \cite{SolovayTennenbaum} works with continuous, commutative systems of complete embeddings which are indexed along a well-order. The following holds: If every forcing of the system 
satisfies a chain condition, then also the direct limit does. Assume for example that all forcings of the system are countable. Then its direct limit satisfies ccc. Assume, moreover, that we want to construct a forcing of size $\omega_2$. Then taking the direct limit will not work, because in our case the limit forcing has size $\leq \omega_1$. To overcome this difficulty, we do not consider a linear system which is indexed along a well-order but a two-dimensional system indexed along a simplified $(\omega _1,1)$-morass. As an example for the approach we constructed a ccc forcing which adds an $\omega_2$-Suslin tree. The conditions of this forcing are Tennenbaum's finite conditions for adding a Suslin tree \cite{Tennenbaum}. However, this forcing does not satisfy ccc on $\omega_2$. Therefore, we apply our approch. That is, our construction uses in every step a countable version of Tennenbaum's forcing, and to obtain complete embeddings we have to thin out these forcings. This results in a thinned out version of Tennenbaum's forcing which satisfies ccc, but still adds an $\omega_2$-Suslin tree.  
\smallskip\\
The kind of two-dimensional system defined in \cite{Irrgang4} is called a FS system along a simplified $(\kappa,1)$-morass. In the present paper, we will generalize the approach to three-dimensional systems, so-called FS systems along simplified $(\kappa,2)$-morasses. We will also observe that under a very weak additional assumption the forcing obtained from a FS system along a simplified gap-1 or gap-2 morass is forcing equivalent to a small subforcing. An immediate consequence of this and \cite{Irrgang4} is: If there is a simplified $(\omega_1,1)$-morass, then there exists a ccc forcing of size $\omega_1$ that adds an $\omega_2$-Suslin tree.  This improves theorem 7.5.1. in Todorcevic's book \cite{StevoBook}: There exists consistently a ccc forcing which adds an $\omega_2$-Suslin tree.    
\smallskip\\
The main result is: If there is a simplified $(\omega_1,2)$-morass, then there exists a ccc forcing of size $\omega_1$ that adds  a $0$-dimensional Hausdorff topology $\tau$ on $\omega_3$ which has spread $s(\tau)=\omega_1$. This forcing is obtained by a FS system along a simplified $(\omega_1,2)$-morass. Its conditions are finite functions $p:x_p\rightarrow 2$ with $x_p \subseteq \omega_3\times \omega_2$. By a theorem of Hajnal and Juhasz \cite{HajnalJuhasz}, $card(X) \leq 2^{2^{s(X)}}=exp(exp(s(X))$ holds for all Hausdorff spaces $X$. In \cite{Juhasz}, Juhasz explicitly raises  the question if the second $exp$ is really necessary. By the usual argument used for Cohen forcing, a ccc forcing of size $\omega_1$ preserves $GCH$. Hence our result shows that it is consistent that there exists a $0$-dimensional Hausdorff space $X$ with $s(X)=\omega_1$ such that $card(X)=2^{2^{s(X)}}$. So far, the consistency of $card(X)=2^{2^{s(X)}}$ has only been known for the case $s(X)=\omega$. The example is the $0$-dimensional, hereditarily separable, hereditarily normal space constructed from $\diamondsuit$ by Fedorcuk \cite{Fedorcuk}.  The author would like to thank Professor Juhasz for pointing this out to him.
\smallskip\\
While the general method of FS systems can be generalized straightforwardly to higher dimensions, we cannot expect that the consistency statements can naively be extended by raising the cardinal parameters. In particular, we cannot expect to be a able to construct from a $(\omega_1,3)$-morass a ccc forcing of size $\omega_1$ which adds a $T_2$ space of size $\omega_4$ and spread $\omega_1$. If this was possible, we could find such a forcing in $L$. However, by the usual argument used for Cohen forcing it preserves $GCH$ which contradicts the theorem of Hajnal and Juhasz. The reason why this generalization does not work is that the gap-3  case yields a four-dimensional construction. Therefore, the finite conditions of our forcing have to fit together appropriately in four directions instead of three and that is impossible. So if and how a statement generalizes to higher-gaps depends heavily on the concrete conditions.   
\smallskip\\
The author started to develop the method of forcing along morasses, because he was interested in solving consistency questions like the following for higher cardinals:
Can there exist a superatomic Boolean algebra with width $\omega$ and height $\omega_2$ (Baumgartner and Shelah \cite{BaumgartnerShelah}, Martinez \cite{Martinez})?
Is it possible that there is a function $f:\omega_2 \times \omega_2 \rightarrow \omega$ such that $f$ is not constant on any rectangle with infinite sides (Todorcevic \cite{Todorcevic,StevoBook})? However, the existence of such a Boolean algebra as well as the existence of such a function contradicts $GCH$. So to get the consistencies we have to destroy $GCH$. Hence a simple application of FS systems will not work because of the properties we described above. Therefore, we will introduce so-called $local$ FS systems along simplified morasses.
\smallskip\\
Local FS systems along morasses are also a step forward into another direction: As outlined above, FS systems have obviously a lot in common with finite support iterations. However, this is not true for all properties of FS iterations. Most prominently,  if $\P$ is the limit of a finite support iteration indexed along $\alpha$, then we can understand a $\P$-generic extension as being obtained successively in $\alpha$-many steps. Moreover, there are names for the forcings used in the single steps. In the case of FS systems, it is unclear what a similar analysis looks like, but if we had it, it would be completely justified to think of our constructions as higher-dimensional FS forcing iterations. 
\smallskip\\
The idea of local FS systems is as follows: Assume that $\langle \P_\eta \mid \eta \leq \kappa^+\rangle$ is a normal, linear FS iteration given as a set of $\kappa^+$-sequences $p \in \P_{\kappa^+}$ such that $\P_\eta=\{ p\upharpoonright \eta \mid p \in \P_{\kappa^+}\}$ and $\P_{\eta+1}\cong\P_\eta \ast \dot Q_\eta$ (where $\dot Q_\eta$ is a $\P_\eta$-name such that $\P_\eta \Vdash(\dot Q_\eta$ is a forcing)). Then $p:\kappa^+ \rightarrow V\in \P_{\kappa^+}$ iff $\P_\eta \Vdash p(\eta) \in \dot Q_\eta$ for all $\eta \in \kappa^+$ and $supp(p):=\{ \eta \in \kappa^+ \mid \P_\eta \not\Vdash p(\eta)=1_{\dot Q_\eta}\}$ is finite. Now, assume that every $\P_\Delta:=\{ p \in \P_{\kappa^+} \mid supp(p) \subseteq \Delta\}$ is obtained through a FS system and therefore satifies a chain condition. Then $\P_{\kappa^+}$ also does. 
\smallskip\\
So far, we do not know how to actually do this with names $\dot Q_\eta$. However, we will give an easy example where no names are needed. Namely, we construct a ccc forcing which adds a chain $\langle X_\alpha \mid \alpha < \omega_2\rangle$ such that $X_\alpha \subseteq \omega_1$, $X_\beta -X_\alpha$ is finite and $X_\alpha-X_\beta$ has size $\omega_1$ for all $\beta < \alpha < \omega_2$. Koszmider constructs such a forcing in \cite{Koszmider} using a Todorcevic $\rho$-function.   
 \smallskip\\
 Todorcevic's method of $\rho$-functions and Shelah's historicized forcing \cite{BaumgartnerShelah,ShelahStanley} seem to be closely related to our approach. Todorcevic uses walks on ordinals to construct $\rho$-functions. A detailed account on the method is his book \cite{StevoBook}. The exact relationship between the two mentioned methods and FS systems is however unclear and would definitely be worth studying. To the author's knowledge, the only result in this direction is by  Morgan \cite{Morgan1996}. He shows that it is possible to directly read off a $\rho$-function from a simplified gap-1 morass. If we use this $\rho$-function and define a forcing to add  a chain $\langle X_\alpha \mid \alpha < \omega_2\rangle$ such that $X_\alpha \subseteq \omega_1$, $X_\beta -X_\alpha$ is finite and $X_\alpha-X_\beta$ has size $\omega_1$ for all $\beta < \alpha < \omega_2$ like Koszmider, then we get exactly the same forcing as with our approach.  
 \smallskip\\
Morasses were introduced by  Jensen in the early 1970's to solve the cardinal transfer problem of model theory in $L$ (see e.g. Devlin \cite{Devlin}). For the proof of the gap-2 transfer theorem a gap-1 morass is used. For higher-gap transfer theorems Jensen has developed so-called higher-gap morasses \cite{Jensen1}. In his Ph.D. thesis, the author generalized these to gaps of arbitrary size \cite{Irrgang3,Irrgang1,Irrgang2}. The theory of morasses is very far developed and very well examined. In particular it is known how to construct morasses in $L$ \cite{Devlin,Friedman,Irrgang3,Irrgang2} and how to force them \cite{Stanley2,Stanley}. Moreover,  Velleman has defined so-called simplified morasses, along which morass constructions can be carried out very easily compared to classical morasses  \cite{Velleman1984,Velleman1987a,Velleman1987b}. Their existence is equivalent to the existence of usual morasses \cite{Donder,Morgan1998}. The fact that the theory of morasses is so far developed is an advantage of the morass approach compared to historic forcing or $\rho$-functions. It allows canonical generalizations to higher cardinals, as shown below. 
\smallskip\\
Finally, we should also mention that besides historicized forcing and $\rho$-functions there is another, quite different method to prove consistencies in two-cardinal combinatorics. This is the method of forcing with models as side conditions or with side conditions in morasses. Models as side conditions were introduced by  Todorcevic \cite{Todorcevic1985,Todorcevic1989}, which was further developed by  Koszmider \cite{Koszmider2000} to side conditions in morasses. Unlike the other methods, it produces proper forcings which are usually not ccc. This is sometimes necessary. For example, Koszmider proved that if CH holds, then there is no ccc forcing that adds a sequence of $\omega_2$ many functions $f:\omega_1 \rightarrow \omega_1$ which is ordered by strict domination mod finite. However, he is able to produce a proper forcing which adds such a sequence \cite{Koszmider2000}. More on the method, including a discussion of its relationship with that of using $\rho$-functions, can be found in Morgan's paper \cite{Morgan2006}. In the context of our approach, this raises the question if it is possible to define something like a countable support iteration along a morass. 
 \section{Simplified gap-2 morasses}
 In this section, we will recall the definition of simplified gap-2 morasses and summarize their properties to the extent necessary for our applications. Except for theorem 2.3 (a) and lemma 2.6 (7), all results in this section are due to Velleman \cite{Velleman1984,Velleman1987a}. Nevertheless, we will usually quote the author's paper \cite{Irrgang4} on FS systems along gap-1 morasses instead of \cite{Velleman1984}, because we hope that in this way the connection to FS systems becomes clearer.
 \smallskip\\
A simplified $(\kappa, 1)$-morass is a structure $\frak{M}=\langle \langle \theta _\alpha \mid \alpha \leq \kappa \rangle , \langle \frak{F}_{\alpha\beta}\mid \alpha < \beta \leq \kappa\rangle\rangle$ satisfying the following conditions:
\smallskip\\
(P0) (a) $\theta _0=1$, $\theta _\kappa =\kappa^+$, $\forall \alpha < \kappa \ \ 0<\theta_\alpha < \kappa$.
\smallskip\\
(b) $\frak{F}_{\alpha\beta}$ is a set of order-preserving functions $f:\theta _\alpha \rightarrow \theta _\beta$.
\smallskip\\
(P1) $|\frak{F}_{\alpha\beta}| < \kappa$ for all $\alpha < \beta < \kappa$.
\smallskip\\
(P2) If $\alpha < \beta < \gamma$, then $\frak{F}_{\alpha\gamma}=\{ f \circ g \mid f \in \frak{F}_{\beta\gamma}, g \in \frak{F}_{\alpha\beta}\}$.
\smallskip\\
(P3) If $\alpha < \kappa$, then $\frak{F}_{\alpha, \alpha +1}=\{ id \upharpoonright \theta _\alpha , h_\alpha\}$ where $h_\alpha$ is such that $h_\alpha \upharpoonright \delta = id \upharpoonright \delta$ and $h_\alpha (\delta) \geq \theta _\alpha$ for some $\delta < \theta _\alpha$.
\smallskip\\
(P4) If $\alpha \leq \kappa$ is a limit ordinal, $\beta _1,\beta_2 < \alpha$ and $f_1 \in \frak{F}_{\beta _1\alpha}$, $f_2 \in \frak{F}_{\beta _2\alpha}$, then there are a $\beta _1, \beta _2 < \gamma < \alpha$, $g \in \frak{F}_{\gamma\alpha}$ and  $j_1 \in \frak{F}_{\beta _1\gamma}$, $j_2 \in \frak{F}_{\beta _2\gamma}$ such that $f_1=g \circ j_1$ and $f_2=g \circ j_2$.
\smallskip\\
(P5) For all $\alpha >0$, $\theta _\alpha =\bigcup \{ f[\theta _\beta] \mid \beta < \alpha , f \in \frak{F}_{\beta\alpha}\}$.
\bigskip\\
Our simplified $(\kappa,1)$-morasses are what are called neat simplified $(\kappa,1)$-morasses in \cite{Velleman1984}. Velleman shows there that if there is one of his simplified $(\kappa,1)$-morasses there is a neat one. Note, moreover, that it is equivalent to replace ``$h_\alpha (\delta) \geq \theta _\alpha$ for some $\delta < \theta _\alpha$'' in (P3) with ``$h_\alpha(\delta+\eta)=\theta_\alpha+\eta$ for some $\delta < \theta_\alpha$ and all $\eta$ such that $\delta+\eta < \theta_\alpha$''. This is easily seen using (P5) and (P2).
\pagebreak\\
{\bf Lemma 2.1}
\smallskip\\
Let $\alpha < \beta < \kappa$, $\tau _1,\tau _2 < \theta _\alpha$, $f_1,f_2 \in \frak{F}_{\alpha \beta}$ and $f_1(\tau _1)=f_2(\tau _2)$. Then $\tau _1=\tau _2$ and $f_1\upharpoonright \tau _1 = f_2 \upharpoonright \tau _2$.
\smallskip\\
{\bf Proof:} See \cite{Irrgang4}, lemma 3.1. $\Box$   
\bigskip\\
A simplified morass defines a tree $\langle T , \prec \rangle$.
\medskip\\
Let $T=\{ \langle \alpha , \nu\rangle \mid \alpha \leq \kappa , \nu < \theta _\alpha\}$.
\smallskip\\
For $t =\langle \alpha , \nu \rangle \in T$ set $\alpha (t)=\alpha$ and $\nu (t)=\nu$.
\smallskip\\
Let $\langle\alpha , \nu \rangle \prec \langle\beta , \tau\rangle$ iff $\alpha <\beta$ and $f(\nu)=\tau$ for some $f \in \frak{F}_{\alpha \beta}$.
\smallskip\\
If $s \prec t$, then $f \upharpoonright (\nu (s) +1)$ does not depend on $f$ by lemma 2.1. So we may define $\pi _{st}:=f \upharpoonright (\nu (s) +1)$.
\bigskip\\
{\bf Lemma 2.2}
\smallskip\\
The following hold:
\smallskip\\
(a) $\prec$ is a tree, $ht_T(t)=\alpha (t)$.
\smallskip\\
(b) If $t_0 \prec t_1 \prec t_2$, then $\pi _{t_0t_1}=\pi _{t_1t_2} \circ \pi _{t_0t_1}$.
\smallskip\\
(c) Let $s \prec t$ and $\pi =\pi _{st}$. If $\pi (\nu^\prime)=\tau^\prime$, $s^\prime=\langle \alpha (s), \nu ^\prime\rangle$ and $t^\prime=\langle \alpha (t), \tau^\prime\rangle$, then $s^\prime \prec t^\prime$ and $\pi _{s^\prime t^\prime}=\pi \upharpoonright (\nu^\prime +1)$.
\smallskip\\
(d) Let $\gamma \leq \kappa$, $\gamma \in Lim$. Let $t \in T_\gamma$. Then $\nu (t) +1=\bigcup \{ rng(\pi_{st})\mid s \prec t\}$.
\smallskip\\
{\bf Proof:} See \cite{Irrgang4}, lemma 3.2. $\Box$
\bigskip\\
A fake gap-1 morass is a structure $\langle \langle \varphi _\zeta \mid \zeta \leq \theta \rangle , \langle \frak{G}_{\zeta\xi}\mid \zeta < \xi \leq \theta \rangle\rangle$ which satisfies the definition of simplified gap-1 morass, except that $\theta$ need not be a cardinal and there is no restriction on the cardinalities of $\varphi _\zeta$ and $\frak{G}_{\zeta\xi}$. Let $\frak{G}_{\zeta,\zeta +1}=\{ id,b\}$. Then the critical point of $b$ is denoted by $\delta_\zeta$ and called the split (or splitting) point of  $\frak{G}_{\zeta,\zeta +1}=\{ id,b\}$. 
\medskip\\
Suppose that  $\langle \langle \varphi _\zeta \mid \zeta \leq \theta \rangle , \langle \frak{G}_{\zeta\xi}\mid \zeta < \xi \leq \theta \rangle\rangle$ and $\langle \langle \varphi^\prime _\zeta \mid \zeta \leq \theta^\prime \rangle , \langle \frak{G}^\prime_{\zeta\xi}\mid \zeta < \xi \leq \theta^\prime \rangle\rangle$ are fake gap-1 morasses. An embedding from the first one to the second will be a function $f$ with domain 
$$(\theta +1) \cup \{ \langle \zeta , \tau \rangle \mid \zeta \leq \theta , \tau < \varphi _\zeta\} \cup \{ \langle \zeta , \xi , b \rangle \mid \zeta < \xi \leq \theta , b \in \frak{G}_{\zeta\xi}\}$$
satisfying certain requirements. We will write $f_\zeta (\tau)$ for $f(\langle \zeta , \tau\rangle)$ and $f_{\zeta\xi}(b)$ for $f(\langle\zeta , \xi ,b\rangle)$.
\medskip\\
The properties are the following ones:
\smallskip\\
(1) $f \upharpoonright (\theta +1)$ is an order preserving function from $\theta +1$ to $\theta ^\prime +1$ such that $f(\theta)=\theta ^\prime$.
\smallskip\\
(2) For all $\zeta \leq \theta$, $f_\zeta$ is an order preserving function from $\varphi _\zeta$ to $\varphi^\prime_{f(\zeta)}$.
\smallskip\\
(3) For all $\zeta < \xi \leq \theta$, $f_{\zeta\xi}$ maps $\frak{G}_{\zeta\xi}$ to  $\frak{G}^\prime_{f(\zeta)f(\xi)}$.
\smallskip\\
(4) If $\zeta < \theta$, then $f_\zeta(\delta _\zeta)=\delta^\prime_{f(\zeta)}$.
\smallskip\\
(5) If $\zeta < \xi \leq \theta$, $b \in \frak{G}_{\zeta\xi}$ and $c \in \frak{G}_{\xi\eta}$, then $f_{\zeta\eta}(c \circ b)=f_{\xi\eta}(c) \circ f_{\zeta\xi}(b)$.
\smallskip\\
(6) If $\zeta < \xi \leq \theta$ and $b \in \frak{G}_{\zeta\xi}$, then $f_\xi \circ b=f_{\zeta\xi}(b) \circ f_\zeta$.
\medskip\\
Assume in the following that $\theta < \theta^\prime$, $\varphi^\prime_\zeta=\varphi_\zeta$ for $\zeta \leq \theta$ and $\frak{G}^\prime_{\zeta\xi}=\frak{G}_{\zeta\xi}$ for $\zeta < \xi \leq \theta$. And let for the moment being $f \upharpoonright \theta =id$, $f_\zeta=id$ for all $\zeta < \theta$ and $f_{\zeta\xi}=id$ for all $\zeta < \xi < \theta$. Let $f_\theta \in \frak{G}^\prime_{\theta\theta^\prime}$. Then we can define an embedding as follows: If $\zeta < \theta$ and $b \in \frak{G}_{\zeta\theta}$, then $f_{\zeta\theta}(b)=f_\theta \circ b$. We call such an embedding $f$ a left-branching embedding. There are many left-branching embeddings, one for every choice of $f_\theta$.
\medskip\\
An embedding $f$ is right-branching if for some $\eta < \theta$,
\smallskip\\
(1) $f \upharpoonright \eta =id$
\smallskip\\
(2) $f(\eta +\zeta)=\theta + \zeta$ if $\eta + \zeta \leq \theta$
\smallskip\\
(3) $f_\zeta = id$ for $\zeta < \eta$
\smallskip\\
(4) $f_{\zeta\xi}=id$ for $\zeta < \xi<\eta$
\smallskip\\
(5) $f_\eta \in \frak{G}_{\eta\theta}$
\smallskip\\
(6) $f_{\zeta\xi}[\frak{G}_{\zeta\xi}]=\frak{G}^\prime_{f(\zeta)f(\xi)}$ if $\eta \leq \zeta < \xi \leq \theta$.
\medskip\\
An amalgamation is a family of embeddings that contains all possible left-branching embeddings, exactly one right-branching embedding and nothing else. The right-branching embedding corresponds to the maps $h_\alpha$ from (P3) in the gap-1 case. Therefore, we will usually denote it by $h$.
\medskip\\
Let $\kappa \geq \omega$ be regular and $\langle \langle \varphi _\zeta \mid \zeta \leq \kappa^+\rangle , \langle \frak{G}_{\zeta\xi} \mid \zeta < \xi \leq \kappa ^+\rangle\rangle$ a simplified $(\kappa^+,1)$-morass such that  $\varphi_\zeta < \kappa$ for all $\zeta < \kappa$. Let $\langle \theta _\alpha \mid \alpha < \kappa \rangle$ be a sequence such that $0<\theta _\alpha<\kappa$ and $\theta _\kappa =\kappa ^+$. Let $\langle \frak{F}_{\alpha\beta}\mid \alpha < \beta \leq \kappa\rangle$ be such that $\frak{F}_{\alpha \beta}$ is a family of embeddings from $\langle \langle \varphi _\zeta \mid \zeta \leq \theta_\alpha \rangle , \langle \frak{G}_{\zeta\xi}\mid \zeta < \xi \leq \theta_\alpha \rangle\rangle$ to $\langle \langle \varphi _\zeta \mid \zeta \leq \theta_\beta \rangle , \langle \frak{G}_{\zeta\xi}\mid \zeta < \xi \leq \theta_\beta \rangle\rangle$.
\medskip\\
This is a simplified $(\kappa,2)$-morass if it has the following properties:
\smallskip\\
(1) $|\frak{F}_{\alpha\beta}| < \kappa$ for all $\alpha < \beta < \kappa$.
\smallskip\\
(2) If $\alpha < \beta < \gamma$, then $\frak{F}_{\alpha\gamma}=\{ f \circ g \mid f \in \frak{F}_{\beta\gamma}, g \in \frak{F}_{\alpha\beta}\}$. Here $f \circ g$ is the composition of the embeddings $f$ and $g$, which are defined in the obvious way: $(f \circ g)_\zeta =f_{g(\zeta)}\circ g_\zeta$ for $\zeta \leq \theta _\alpha$ and $(f \circ g)_{\zeta\xi}=f_{g(\zeta)g(\xi)}\circ g_{\zeta\xi}$ for $\zeta < \xi \leq \theta _\alpha$.
\smallskip\\
(3) If $\alpha < \kappa$, then $\frak{F}_{\alpha, \alpha +1}$ is an amalgamation.
\smallskip\\
(4) If $\alpha \leq \kappa$ is a limit ordinal, $\beta _1,\beta_2 < \alpha$ and $f_1 \in \frak{F}_{\beta _1\alpha}$, $f_2 \in \frak{F}_{\beta _2\alpha}$, then there are a $\beta _1, \beta _2 < \gamma < \alpha$, $g \in \frak{F}_{\gamma\alpha}$ and  $j_1 \in \frak{F}_{\beta _1\gamma}$, $j_2 \in \frak{F}_{\beta _2\gamma}$ such that $f_1=g \circ j_1$ and $f_2=g \circ j_2$.
\smallskip\\
(5) For all $\alpha \leq \kappa$, $\alpha \in Lim$:
\smallskip\\
(a) $\theta _\alpha =\bigcup \{ f[\theta _\beta] \mid \beta < \alpha , f \in \frak{F}_{\beta\alpha}\}$.
\smallskip\\
(b) For all $\zeta \leq \theta_\alpha$, $\varphi_\zeta =\bigcup \{ f_{\bar\zeta}[\varphi_{\bar\zeta}] \mid \exists \beta < \alpha \ (f \in \frak{F}_{\beta\alpha}$ and $f(\bar\zeta)=\zeta)\}$.
\smallskip\\
(c) For all $\zeta < \xi \leq \theta _\alpha$, $\frak{G}_{\zeta\xi}=\bigcup \{ f_{\bar\zeta\bar\xi}[\frak{G}_{\bar\zeta\bar\xi}]\mid \exists \beta < \alpha \ (f \in \frak{F}_{\beta\alpha}, f(\bar\zeta)=\zeta$ and $ f(\bar\xi)=\xi)\}$.
\bigskip\\
{\bf Theorem 2.3}
\smallskip\\
(a) If $V=L$, then there is a simplified $(\kappa ,2)$-morass for all regular $\kappa >\omega$.
\smallskip\\
(b) If $\kappa >\omega$ is regular, then there is a forcing $\P$ which preserves cardinals and cofinalities such $\P \Vdash$ (there is a simplified $(\kappa , 2)$-morass). 
\smallskip\\
{\bf Proof:} (a) The existence of a gap-2 morass was first proved by Jensen. The proof is very similiar to the existence proof for gap-1 morasses. See Devlin \cite{Devlin}, VIII 2. A sketch of the proof can be found in Friedman \cite{Friedman}, 1.3. That a simplified gap-2 morass can be obtained from an ordinary one was shown by Morgan in \cite{Morgan1998}. 
\smallskip\\
(b) See Velleman  \cite{Velleman1987a}.  $\Box$
\bigskip\\
Since $\langle \langle \varphi _\zeta \mid \zeta \leq \kappa^+\rangle , \langle \frak{G}_{\zeta\xi} \mid \zeta < \xi \leq \kappa ^+\rangle\rangle$ is a simplified $(\kappa^+,1)$-morass, there is a tree $\langle T,\prec\rangle$ with levels $T_\eta$ for $\eta \leq \kappa ^+$ as in lemma 1.2. And there are maps $\pi _{st}$ for $s \prec t$. Moreover, if we set $\frak{F}^\prime_{\alpha\beta}=\{ f \upharpoonright \theta _\alpha \mid f \in \frak{F}_{\alpha\beta}\}$, then $\langle \langle \theta _\alpha \mid \alpha \leq \kappa \rangle , \langle \frak{F}^\prime_{\alpha\beta} \mid \alpha < \beta \leq \kappa \rangle \rangle$ is a simplified $(\kappa ,1)$-morass. So there is also a tree $\langle T^\prime ,\prec^\prime\rangle$ with levels $T^\prime _\eta$ for $\eta \leq \kappa$ as in lemma 2.2 on this morass. Improving lemma 2.1, the following holds:
\bigskip\\
{\bf Lemma 2.4}
\smallskip\\
Suppose $\alpha < \beta \leq \kappa$, $f_1,f_2 \in \frak{F}_{\alpha\beta}$, $\zeta _1,\zeta _2 < \theta _\alpha$ and $f_1(\zeta_1)=f_2(\zeta_2)$. Then $\zeta _1=\zeta _2$, $f_1 \upharpoonright \zeta _1 =f_2 \upharpoonright \zeta _1$, $(f_1)_\xi =(f_2)_\xi$ for all $\xi \leq \zeta_1$, and $(f_1)_{\xi\eta}=(f_2)_{\xi\eta}$ for all $\xi < \eta \leq \zeta _1$.
\smallskip\\
{\bf Proof:} See Velleman \cite{Velleman1987a}, lemma 2.2. $\Box$ 
\bigskip\\
Now, let $s=\langle \alpha , \nu \rangle \in T^\prime _\alpha$, $t=\langle \beta , \tau \rangle \in T^\prime _\beta$ and $s \prec ^\prime t$. Then there is some $f \in \frak{F}^\prime_{\alpha \beta}$ such that $f(\nu)=\tau$. By lemma 2.4 
$$f \upharpoonright ((\nu +1) \cup \{ \langle \zeta , \tau \rangle \mid \zeta \leq \nu , \tau < \varphi _\zeta\} \cup \{ \langle \zeta , \xi , b \rangle \mid \zeta < \xi \leq \nu , b \in \frak{G}_{\zeta\xi}\} )$$
does not depend on $f$. So we may call it $\pi^\prime _{st}$.
\bigskip\\
Finally, we can prove something very natural:
\medskip\\
{\bf Lemma 2.5}
\smallskip\\
(a) If $\zeta < \xi \leq \kappa^+$, then $id \upharpoonright \varphi _\zeta \in \frak{G}_{\zeta\xi}$.
\smallskip\\
(b) If $\alpha < \beta \leq \kappa$, then there is a $g \in \frak{F}_{\alpha\beta}$ such that $g \upharpoonright \theta _\alpha = id \upharpoonright \theta _\alpha$. \smallskip\\
{\bf Proof:} (a) See \cite{Irrgang4}, lemma 3.3. 
\smallskip\\
(b) See Velleman  \cite{Velleman1987a}, lemma 2.4. $\Box$
\bigskip\\
In addition to the maps $f \in \frak{F}_{\alpha\beta}$, we need maps $\bar f$ that are associated to $f$. For a set of ordinals $X$, let $ssup(X)$ be the least $\alpha$ such that $X \subseteq \alpha$. And let $\bar f(\zeta)=ssup(f[\zeta]) \leq f(\zeta)$.
\medskip\\
{\bf Lemma 2.6}
\smallskip\\
For every $\alpha < \beta \leq \kappa$, $f \in \frak{F}_{\alpha\beta}$ and $\zeta \leq \theta_\alpha$, there are unique functions $\bar f_\zeta : \varphi _\zeta \rightarrow \varphi _{\bar f(\zeta)}$, $\bar f_{\xi\zeta}:\frak{G}_{\xi\zeta} \rightarrow \frak{G}_{f(\xi)\bar f(\zeta)}$ for all $\xi < \zeta$, and $f^\#(\zeta) \in \frak{G}_{\bar f(\zeta)f(\zeta)}$ such that:
\smallskip\\
(1) $f_\zeta=f^\#(\zeta) \circ \bar f_\zeta$
\smallskip\\
(2) $\forall \xi < \zeta \ \forall b \in \frak{G}_{\xi\zeta} \ \ f_{\xi\zeta}(b)=f^\#(\zeta)\circ \bar f_{\xi\zeta}(b)$.
\smallskip\\
Furthermore, these functions have the following properties:
\smallskip\\
(3) If $\xi < \bar f(\zeta)$ and $b \in \frak{G}_{\xi\bar f(\zeta)}$, then $\exists \eta < \zeta \ \exists c \in \frak{G}_{\eta\zeta} \ \exists d \in \frak{G}_{\xi f(\eta)}$ $$b=\bar f_{\eta\zeta}(c) \circ d.$$
(4) $\forall \xi < \zeta \ \forall  b \in \frak{G}_{\xi \zeta}\ \ \bar f_\zeta \circ b =\bar f_{\xi\zeta}(b) \circ f_\xi$.
\smallskip\\
(5) If $\eta < \xi < \zeta$, $b \in \frak{G}_{\xi\zeta}$ and $c \in \frak{G}_{\eta\xi}$, then $\bar f_{\eta\zeta}(b \circ c)=\bar f_{\xi\zeta}(b) \circ f_{\eta\xi}(c)$.
\smallskip\\
(6) If $\alpha < \beta < \gamma \leq \kappa$, $f \in \frak{F}_{\beta\gamma}$, $g \in \frak{F}_{\alpha\beta}$ and $\zeta \leq \theta_\alpha$, then

$(\overline{f \circ g})_\zeta =\bar f_{\bar g(\zeta)} \circ \bar g_\zeta$

$(f \circ g)^\#(\zeta)=f_{\bar g(\zeta)g(\zeta)}(g^\#(\zeta)) \circ f^\#(\bar g(\zeta))$ and

$(\overline{f\circ g})_{\xi\zeta}=\bar f_{g(\xi)\bar g(\zeta)} \circ \bar g_{\xi\zeta}$ for all $\xi < \zeta$.
\smallskip\\
{\bf Proof:}  See Velleman  \cite{Velleman1987a}, lemma 2.1.  $\Box$ 
\bigskip\\
From the previous lemma, we get of course also maps $(\overline{\pi ^\prime_{st}})_\zeta$ for $s \prec ^\prime t$ and $\zeta \leq \nu(t)$.    
\section{FS systems along morasses}
In this section, we recall the definition of FS systems along gap-1 morasses given in \cite{Irrgang4} and generalize it to the gap-2 case, which is straightforward. 
\smallskip\\
Let $\P$ and $\Q$ be partial orders. A map $\sigma :\P \rightarrow \Q$ is called a complete embedding if
\smallskip\\
(1) $\forall  p ,p^\prime \in \P \ (p^\prime \leq p \rightarrow \sigma (p^\prime ) \leq \sigma (p))$
\smallskip\\
(2) $\forall  p ,p^\prime \in \P \ ( p$ and $p^\prime$ are incompatible $\leftrightarrow$ $\sigma (p)$ and $\sigma (p^\prime )$ are incompatible)
\smallskip\\
(3) $\forall q \in \Q \ \exists p \in \P \ \forall p^\prime \in \P \ (p^\prime \leq p \rightarrow (\sigma (p^\prime )$ and $q$ are compatible in $\Q$)).
\smallskip\\
In (3), we call $p$ a reduction of $q$ to $\P$ with respect to $\sigma$.
\medskip\\
If only (1) and (2) hold, we say that $\sigma$ is an embedding. If $\P \subseteq \Q$ such that the identity is an embedding, then we write $\P \subseteq _\bot \Q$.
\medskip\\
We say that $\P \subseteq \Q$ is completely contained in $\Q$ if $id \upharpoonright \P :\P \rightarrow \Q$ is a complete embedding.
\bigskip\\
Let $\langle \langle \varphi _\zeta \mid \zeta \leq \kappa^+ \rangle , \langle \frak{G}_{\zeta\xi}\mid \zeta < \xi \leq \kappa^+ \rangle\rangle$ be a simplified $(\kappa^+ ,1)$-morass. We want to "iterate" along it. This leads to the following definition.
\medskip\\
We say that $\langle\langle \P _\eta \mid \eta \leq \kappa ^{++} \rangle ,\langle \sigma _{st} \mid s \prec t \rangle , \langle e_\alpha \mid \alpha < \kappa^+ \rangle\rangle$ is a FS system along  $\langle \langle \varphi _\zeta \mid \zeta \leq \kappa^+ \rangle , \langle \frak{G}_{\zeta\xi}\mid \zeta < \xi \leq \kappa^+ \rangle\rangle$ if the following conditions hold:
\medskip\\
(FS1) $\langle \P _\eta \mid \eta \leq \kappa ^{++}  \rangle$ is a sequence of partial orders such that $\P _\eta \subseteq _\bot \P _\nu$ if $\eta \leq \nu$ and $\P _\lambda =\bigcup \{ \P _\eta\mid \eta < \lambda\}$ for $\lambda \in Lim$.\medskip\\
(FS2) $\langle \sigma _{st} \mid s \prec t \rangle$ is a commutative system of injective embeddings $\sigma _{st}:\P _{\nu (s)+1} \rightarrow \P _{\nu (t)+1}$ such that if $t$ is a limit point in $\prec$, then 
$$\P _{\nu (t) +1} = \bigcup \{ \sigma _{st}[\P _{\nu (s)+1}]\mid s \prec t \}.$$ 
\smallskip\\
(FS3) $e_\alpha : \P _{\varphi _{\alpha +1}} \rightarrow \P _{\varphi _\alpha}$.
\medskip\\
(FS4) Let $s \prec t$ and $\pi =\pi _{st}$. If $\pi (\nu^\prime)=\tau^\prime$, $s^\prime=\langle \alpha (s), \nu ^\prime\rangle$ and $t^\prime=\langle \alpha (t), \tau^\prime\rangle$, then $\sigma _{st}:\P_{\nu (s)+1} \rightarrow \P _{\nu (t) +1}$ extends $\sigma _{s^\prime t^\prime}:\P _{\nu ^\prime +1}\rightarrow \P _{\tau^\prime +1}$.
\medskip\\
Hence for $f \in \frak{G}_{\alpha\beta}$, we may define $\sigma _f  =\bigcup \{ \sigma _{st} \mid s=\langle \alpha , \nu \rangle, t=\langle \beta , f(\nu)\rangle\}$.
\medskip\\
(FS5) If $\pi _{st}=id \upharpoonright \nu (s)+1$, then $\sigma _{st}=id \upharpoonright \P _{\nu (s) +1}$.
\medskip\\
(FS6)(a) If $\alpha < \kappa^+$, then $\P _{\varphi_\alpha}$ is completely contained in $\P _{\varphi _{\alpha +1}}$ in such a way that $e_\alpha(p)$ is a reduction of $p \in \P _{\varphi _{\alpha +1}}$. 
\smallskip\\
(b) If $\alpha < \kappa^+ $, then $\sigma _\alpha :=\sigma _{h_\alpha} :\P _{\varphi _\alpha} \rightarrow \P _{\varphi _{\alpha +1}}$ is a complete embedding such that $e_\alpha(p)$ is a reduction of $p \in \P _{\varphi _{\alpha +1}}$.
\medskip\\
(FS7)(a) If $\alpha < \kappa^+$ and $p \in \P _{\varphi _\alpha}$, then $e_\alpha(p)=p$.
\smallskip\\
(b) If $\alpha < \kappa^+$ and $p \in rng(\sigma _\alpha)$, then $e_\alpha(p)=\sigma ^{-1}_\alpha(p)$.
\bigskip\\
The definition of an FS system along a simplified $(\kappa,1)$-morass, of course, makes sense for arbitrary regular $\kappa \geq \omega$. We gave it here for successor cardinals because if a simplified $(\kappa,2)$-morass is given then the associated gap-1 morass $\langle \langle \varphi _\zeta \mid \zeta \leq \kappa^+ \rangle , \langle \frak{G}_{\zeta\xi}\mid \zeta < \xi \leq \kappa^+ \rangle\rangle$ is a simplified $(\kappa^+,1)$-morass. 
\bigskip\\
To simplify notation, set $\P := \P _{\kappa ^{++}}$.
\medskip\\
As in the case of (linear) FS iterations it is sometimes more convenient to represent $\P$ as a set of functions $p^\ast:\kappa ^+ \rightarrow V$ such that $p^\ast(\alpha) \in \P _{\varphi _\alpha}$ for all $\alpha < \kappa ^+$.
\bigskip\\
To define such a function $p^\ast$ from $p \in \P$ set recursively

$p_0=p$

$\nu _n(p)=min\{ \eta \mid p_n \in \P _{\eta +1}\}$

$t_n(p)=\langle \kappa^+ ,\nu _n(p)\rangle$

$p^{(n)}(\alpha)=\sigma ^{-1}_{st}(p_n)$ if $s \in T_\alpha$, $s \prec t:= t_n(p)$ and $p_n \in rng(\sigma_{st})$.

Note that, by lemma 2.2 (a), $s$ is uniquely determined by $\alpha$ and $t_n(p)$. Hence we really define a function. Set

$\gamma _n(p)=min(dom(p^{(n)}))$.

By (FS2), $\gamma _n(p)$ is a successor ordinal or $0$. Hence, if $\gamma_n(p)\neq 0$, we may define

$p_{n+1}=e_{\gamma_n(p)-1}(p^{(n)}(\gamma _n(p)))$.

If $\gamma_n(p)= 0$, we let $p_{n+1}$ be undefined.
\medskip\\ 
Finally, set $p^\ast=\bigcup\{ p^{(n)} \upharpoonright [\gamma _n(p) , \gamma _{n-1}(p)[ \ \mid n \in \omega \}$ where $\gamma _{-1}(p)=\kappa^+$.
\medskip\\
Note: If $n>0$ and $\alpha \in [\gamma _n(p),\gamma_{n-1}(p)[$, then $p^\ast(\alpha)=\sigma _{s\bar t}^{-1}(p_n)$ where $\bar t=\langle \gamma_{n-1}(p)-1,\nu _n(p)\rangle$ because $p^\ast(\alpha)=p^{(n)}(\alpha)=\sigma^{-1}_{st}(p_n)=(\sigma _{\bar tt}\circ \sigma _{s\bar t})^{-1}(p_n)=\sigma_{s\bar t}(p_n)$ where $t=t_n(p)=\langle \kappa,\nu_n(p)\rangle$. The first two equalities are just the definitions of $p^\ast$ and $p^{(n)}$. For the third equality note that $\bar t \prec t$ by lemma 2.5 (a). So the equality follows from the commutativity of $\langle \sigma _{st} \mid s \prec t \rangle$. The last equality holds by (FS5).
\medskip\\
It follows from the previous observation that $\langle \gamma _n(p) \mid n \in \omega \rangle$ is decreasing. So the recursive definition above breaks down at some point, i.e. $\gamma _n(p)=0$ for some $n \in \omega$. Hence
$$supp(p)=\{ \gamma _n(p) \mid n \in \omega\}$$
is finite.  
\bigskip\\
{\bf Lemma 3.1}
\smallskip\\
If $p^\ast(\alpha)$ and $q^\ast(\alpha)$ are compatible for $\alpha=max(supp(p) \cap supp(q))$, then $p$ and $q$ are compatible.
\smallskip\\
{\bf Proof:} See \cite{Irrgang4}, lemma 4.1. $\Box$ 
\bigskip\\
{\bf Theorem 3.2}
\smallskip\\
Let $\mu ,\kappa >\omega$ be cardinals, $\kappa$ regular. Let  $\langle\langle \P _\eta \mid \eta \leq \kappa ^+ \rangle ,\langle \sigma _{st} \mid s \prec t \rangle ,\langle e_\alpha \mid \alpha < \kappa \rangle \rangle$ be a FS system along a $(\kappa ,1)$-morass $\frak{M}$. Assume that all $\P _\eta$ with $\eta < \kappa$ satisfy the $\mu$-cc. Then $\P _{\kappa ^+}$ also does.
\smallskip\\
{\bf Proof:} See \cite{Irrgang4}, lemma 4.2. $\Box$
\bigskip\\
Now, let $\frak{M}$ be a simplified $(\kappa,2)$-morass.
\medskip\\
We say that $$\langle\langle \P _\eta \mid \eta \leq \kappa ^{++} \rangle ,\langle \sigma _{st} \mid s \prec t \rangle ,\langle \sigma^\prime _{st} \mid s \prec^\prime t \rangle , \langle e_\alpha \mid \alpha < \kappa^+ \rangle, \langle e^\prime_\alpha \mid \alpha < \kappa \rangle\rangle$$ is a FS system along $\frak{M}$ if the following conditions hold:
\medskip\\
(FS$_2$1) $\langle\langle \P _\eta \mid \eta \leq \kappa ^{++} \rangle ,\langle \sigma _{st} \mid s \prec t \rangle , \langle e_\alpha \mid \alpha < \kappa^+ \rangle\rangle$ is a FS system along  $\langle \langle \varphi _\zeta \mid \zeta \leq \kappa^+ \rangle , \langle \frak{G}_{\zeta\xi}\mid \zeta < \xi \leq \kappa ^+ \rangle\rangle$.
\medskip\\
Let $\Q =\{ p^\ast \upharpoonright supp(p) \mid p \in \P_{\kappa^{++}}\}$.
\smallskip\\
Define a partial order $\leq$ on $\Q$ by setting $p \leq q$ iff $dom(q) \subseteq dom(p)$ and $p (\alpha) \leq q(\alpha)$ for all $\alpha \in dom(q)$.
\smallskip\\
Set $\Q_\gamma := \{ p \in \Q \mid dom(p)\subseteq \gamma\}$.
\medskip\\
(FS$_2$2) $\langle \sigma^\prime _{st} \mid s \prec^\prime t \rangle$ is a commutative system of injective embeddings $\sigma^\prime _{st}:\Q _{\nu (s)+1} \rightarrow \Q _{\nu (t)+1}$ such that if $t$ is a limit point in $\prec ^\prime$, then $\Q _{\nu (t) +1} = \bigcup \{ \sigma^\prime _{st}[\Q _{\nu (s)+1}]\mid s \prec^\prime t \}$. 
\smallskip\\
(FS$_2$3) $e^\prime_\alpha : \Q _{\theta _{\alpha +1}} \rightarrow \Q _{\theta _\alpha}$.
\medskip\\
(FS$_2$4) Let $s \prec^\prime t$ and $\pi =\pi^\prime _{st}$. If $\pi (\nu^\prime)=\tau^\prime$, $s^\prime=\langle \alpha (s), \nu ^\prime\rangle$ and $t^\prime=\langle \alpha (t), \tau^\prime\rangle$, then $\sigma^\prime _{st}:\Q_{\nu (s)+1} \rightarrow \Q _{\nu (t) +1}$ extends $\sigma^\prime _{s^\prime t^\prime}:\Q _{\nu ^\prime +1}\rightarrow \Q _{\tau^\prime +1}$.
\medskip\\
Hence for $f \in \frak{F}_{\alpha\beta}$, we may define $\sigma _f  =\bigcup \{ \sigma _{st} \mid s=\langle \alpha , \nu \rangle, t=\langle \beta , f(\nu)\rangle\}$.
\medskip\\
(FS$_2$5) If $\pi^\prime _{st}\upharpoonright \nu (s)+1=id \upharpoonright \nu (s)+1$, then $\sigma^\prime _{st}=id \upharpoonright \Q _{\nu (s) +1}$.
\medskip\\
(FS$_2$6)(a) If $\alpha < \kappa$, then $\Q _{\theta_\alpha}$ is completely contained in $\Q _{\theta _{\alpha +1}}$ in such a way that $e^\prime_\alpha(p)$ is a reduction of $p \in \Q _{\theta _{\alpha +1}}$. 
\smallskip\\
(b) If $\alpha < \kappa$, then $\sigma^\prime _\alpha :=\sigma^\prime _{h_\alpha} :\Q _{\theta _\alpha} \rightarrow \Q _{\theta _{\alpha +1}}$ (where $h_\alpha$ is the unique right-branching $f \in \frak{F}_{\alpha , \alpha +1}$) is a complete embedding such that $e^\prime_\alpha(p)$ is a reduction of $p \in \Q _{\theta _{\alpha +1}}$.
\medskip\\
(FS$_2$7)(a) If $\alpha < \kappa$ and $p \in \Q _{\theta _\alpha}$, then $e^\prime_\alpha(p)=p$.
\smallskip\\
(b) If $\alpha < \kappa$ and $p \in rng(\sigma^\prime _\alpha)$, then $e^\prime_\alpha(p)=(\sigma^\prime) ^{-1}_\alpha(p)$.
\bigskip\\
This definition deserves some explanation. An FS system along a gap-1 morass is obtained by thinning out a forcing $P$ by recursion along the morass, which yields a forcing $\P$. An example of such a construction is given in \cite{Irrgang4}. Similarly, an FS system along a gap-2 morass is obtained by thinning out a forcing $P$ twice. In the first step, it is thinned out along the gap-1 morass  $\langle \langle \varphi _\zeta \mid \zeta \leq \kappa^+\rangle , \langle \frak{G}_{\zeta\xi} \mid \zeta < \xi \leq \kappa ^+\rangle\rangle$, which yields a forcing $P^\prime$ and an FS system along the gap-1 morass. So it makes sense to consider $Q^\prime = \{ p^\ast \upharpoonright supp(p) \mid p \in P^\prime\}$. Then, in the second step, $P^\prime$ is thinned out to $\P$. This is actually  done by thinning out $Q^\prime$ to the $\Q$ of the definition. This explains why the auxiliary structure is necessary. 
\bigskip\\
{\bf Theorem 3.3}
\smallskip\\
Let $\kappa, \nu > \omega$ be cardinals, $\kappa$ regular. Let $\langle\langle \P _\eta \mid \eta \leq \kappa ^{++} \rangle ,\langle \sigma _{st} \mid s \prec t \rangle ,\langle \sigma^\prime _{st} \mid s \prec^\prime t \rangle , \langle e_\alpha \mid \alpha < \kappa^+ \rangle, \langle e_\alpha \mid \alpha < \kappa \rangle\rangle$ be a FS system along a $(\kappa, 2)$-morass.
\smallskip\\
(a) If $\langle \Q , \leq \rangle$ satisfies the $\mu$-cc, then $\P$ also does.
\smallskip\\
(b) If all $\Q _\eta$ with $\eta < \kappa$ satisfy the $\mu$-cc, then $\P$ also does.
\smallskip\\
{\bf Proof:} (a) follows directly from theorem 3.2.
\smallskip\\
(b) By properties (FS$_2$1) - (FS$_2$7), we obtain as in theorem 3.2, that $\Q$ satisfies the $\mu$-cc. Hence the claim follows by (a). $\Box$
\bigskip\\
If we define $i: \P \rightarrow \Q, p \mapsto p^\ast \upharpoonright supp(p)$ and assume that
\smallskip\\
(1) $\forall p,q \in \P \  \forall \alpha \in \kappa: \ p \leq q \rightarrow e_\alpha(p) \leq e_\alpha(q)$
\smallskip\\
(2) $\forall p \leq q \in \P \ \forall s \prec t: \ p \in rng(\sigma_{st}) \rightarrow q \in rng(\sigma_{st})$,
\smallskip\\
then $i:\P \rightarrow \Q$ is a dense embedding, i.e. $\P$ and $\Q$ are forcing equivalent. Hence in this case, (a) is trivial. As we will see, this is also the reason why the method can hardly be used to construct forcings which destroy $GCH$.
\section{Cohen forcing and a topological space}
To understand how FS systems along morasses work, we will discuss the simplest example, Cohen forcing. That is, we consider the forcing
$$P=\{p:x_p \rightarrow 2 \mid x_p \subseteq \omega_3 \times \omega_2 \hbox{ finite}\}.$$
As usual, we set $p \leq q$ iff $q \subseteq p$.
\smallskip\\
"Iterating" Cohen forcing along a gap-2 morass as in the definion of FS system, will yield  a ccc forcing of size $\omega_1$ that adds a $0$-dimensional $T_2$ topology on $\omega_3$ with spread $\omega_1$. The construction has two important precursors. Those are, firstly, the construction of a ccc forcing that adds an $\omega_2$-Suslin tree in \cite{Irrgang4} and, secondly, Velleman's proof \cite{Velleman1987a} that the model theoretic gap-3 theorem holds in $L$. In the following, we will refer to \cite{Irrgang4} and \cite{Velleman1987a} from time to time to point out similarities between the constructions.  We hope that this makes the whole proof more comprehensible. 
\medskip\\
Let $\pi :\bar \theta \rightarrow \theta$ be an order-preserving map. Then $\pi :\bar\theta \rightarrow \theta$ induces maps $\pi : \bar\theta \times \omega_2 \rightarrow \theta \times\omega_2$ and  $\pi :(\bar \theta  \times \omega_2)\times 2  \rightarrow (\theta  \times \omega_2 ) \times 2$ in the obvious way:
$$\pi : \bar\theta \times \omega_2 \rightarrow \theta \times \omega_2,\quad \langle \gamma ,\delta\rangle \mapsto \langle \pi (\gamma ),\delta \rangle$$
$$\pi :(\bar \theta  \times \omega_2) \times 2 \rightarrow (\theta  \times \omega_2 ) \times 2, \quad \langle x ,\epsilon\rangle \mapsto \langle \pi (x ),\epsilon \rangle.$$
Basically, we will define the maps $\sigma$ of the FS system by setting $\sigma(p)=\pi[p]$.
\bigskip\\
Now, we start our construction of $\P$. In the first step, we define partial orders $P(\tau)$ for $\tau \leq \omega _3$ and $Q(\tau)$ for $\tau \leq \omega _2$. In the second step, we thin out $P(\tau)$ and $Q(\tau)$ to the $\P_\tau$ and $\Q_\tau$ which form the FS system along the gap-2 morass. 
\medskip\\
Assume that a simplified $(\omega_1,2)$-morass as in the previous section is given. We define $P(\tau)$ by induction on the levels of $\langle \langle \varphi_\zeta \mid \zeta \leq \omega _2\rangle,\langle \frak{G}_{\zeta\xi} \mid \zeta < \xi \leq \omega _2\rangle \rangle$ which we enumerate by $\beta \leq \omega_2$. 
\medskip\\
{\it Base Case}: $\beta =0$
\medskip\\
Then we only need to define $P(1)$.
\smallskip\\
Let $P(1):=\{ p \in P \mid x_p \subseteq 1\times \omega\}$.  
\medskip\\
{\it Successor Case}: $\beta = \alpha +1$
\medskip\\
We first define $P(\varphi_\beta)$. Let it be the set of all $p \in P$ such that
\smallskip\\
(1) $x_p \subseteq \varphi_\beta \times \omega\beta$
\smallskip\\
(2) $p \upharpoonright (\varphi_\alpha \times \omega\alpha)$, $h_\alpha^{-1}[p \upharpoonright (\varphi_\beta \times \omega\alpha)]$ $\in P(\varphi_\alpha)$
\smallskip\\
(3) $p \upharpoonright (\varphi_\alpha \times \omega\alpha)$ and $h_\alpha^{-1}[p \upharpoonright (\varphi_\beta \times \omega\alpha)]$ are compatible in $P$
\smallskip\\
where $h_\alpha$ is as in (P3) in the definition of a simplified gap-1 morass. 
\medskip\\
For all $\nu \leq \varphi_\alpha$ $P(\nu)$ is already defined. For $\varphi_\alpha <\nu \leq \varphi_\beta$ set 
$$P (\nu )=\{ p \in P(\varphi_\beta) \mid x_p \subseteq \nu  \times \omega\beta \}.$$
Set 
$$\sigma _{st}:P(\nu (s)+1) \rightarrow P(\nu(t)+1) , p \mapsto \pi _{st}[p].$$ 
It remains to define $e_\alpha$. If $p \in rng(\sigma _\alpha)$, then set $e_\alpha (p)=\sigma^{-1}_\alpha(p)$. If $p \in P (\varphi _\alpha)$, then set $e_\alpha(p)=p$. And if $p \notin rng(\sigma_\alpha) \cup P(\varphi_\alpha)$, then set  
$$e_\alpha(p)= p \upharpoonright (\varphi_\alpha \times \omega\alpha) \cup h_\alpha^{-1}[p \upharpoonright (\varphi_\beta \times \omega\alpha)].$$  
{\it Limit Case}: $\beta \in Lim$
\medskip\\
For $t \in T_\beta$ set $P(\nu(t)+1)=\bigcup \{ \sigma _{st}[P(\nu (s)+1)] \mid s \prec t \}$ and $P (\lambda) =\bigcup \{ P (\eta) \mid \eta < \lambda\}$ for $\lambda \in Lim$ where $\sigma _{st}:P (\nu (s)+1) \rightarrow P(\nu(t)+1), p \mapsto \pi _{st}[p]$.
\bigskip\\
{\bf Lemma 4.1}
\smallskip\\
 $\langle\langle P (\eta) \mid \eta \leq \omega _3 \rangle ,\langle \sigma _{st} \mid s \prec t \rangle , \langle e_\alpha \mid \alpha < \omega_2 \rangle\rangle$ is a FS system along  $\langle \langle \varphi _\zeta \mid \zeta \leq \omega_2 \rangle , \langle \frak{G}_{\zeta\xi}\mid \zeta < \xi \leq \omega_2 \rangle\rangle$.
\smallskip\\
{\bf Proof:} Most things are clear. We only prove (FS6). Let $p \in P(\varphi _\beta)$ and $\beta = \alpha +1$. Let $q:= p \upharpoonright (\varphi_\alpha \times \alpha) \cup h_\alpha^{-1}[p \upharpoonright (\varphi_\beta \times \alpha)]$. We have to prove that $q$ is a reduction of $p$ with respect to $\sigma _\alpha$ and $id \upharpoonright P (\varphi _\alpha)$. To do so, let $r \leq q$. We have to find an $s \leq p,\sigma _\alpha(r),r$ such that $s \in P (\varphi _\beta)$. Define $s$ as $s:=p \cup r \cup h_\alpha[r]$. It is easily seen that $s$ is as wanted. $\Box$ 
\bigskip\\
By the previous lemma every $p \in P(\omega _3)$ has finite support and we may define
\smallskip\\
$p^\ast$ for $p \in P(\omega _3)$ as in section 3
\smallskip\\
$Q=\{ p^\ast \upharpoonright supp(p) \mid p \in P(\omega _3)\}$
\smallskip\\
$Q(\gamma)=\{ p \in Q \mid dom(p) \subseteq \gamma\}$.
\bigskip\\
{\bf Lemma 4.2}
\smallskip\\
If $p \leq q$ in $P(\omega_3)$, then $p^\ast \upharpoonright supp(p)\leq q^\ast \upharpoonright supp(q)$ in $Q$.
\smallskip\\
{\bf Proof:} Let $\nu_0(q) \leq \nu_0(p)$ and $\gamma_0(p)$ be as in the definition of the support of a condition. Let $s \prec t :=t_0(p)$, $s \in T_{\gamma_0(p)}$ and $s^\prime \prec t^\prime:=t_0(q)$, $s^\prime \in T_{\gamma_0(p)}$. Then $\nu_0(q)\in rng(\pi_{st})$ and $\pi_{s^\prime t^\prime}=\pi_{st}(\nu(s^\prime+1)$ by lemma 2.2 (c). Hence $p^\ast (\gamma_0(p)) \leq q^\ast (\gamma_0(p))$ and $\alpha \notin supp(q)$ for all $\gamma_0(p) < \alpha < \omega_2=\gamma_{-1}(p)$. From  $p^\ast (\gamma_0(p)) \leq q^\ast (\gamma_0(p))$ it follows that  $p^\ast (\gamma_0(p)-1) \leq q^\ast (\gamma_0(p)-1)$ by the definition of $\langle e_\alpha \mid \alpha \in \omega_2\rangle$. Now we can repeat this argumentation finitely many times which yields that $supp(q) \subseteq supp(p)$ and that  $p^\ast (\gamma_n(p)) \leq q^\ast (\gamma_n(p))$ for all $n \in \omega$ where it is defined. Hence $p^\ast \upharpoonright supp(p)\leq q^\ast \upharpoonright supp(q)$ as wanted. $\Box$ 
\bigskip\\
So far, the development is as in \cite{Irrgang4}. Following the definitions of section 3, we have to do the same for $Q$. We will, however, not use the maps $f \in \frak{F}_{\alpha\beta}$ but $\bar f$ to map $p \in Q(\theta_\alpha)$ to $Q(\theta_\beta)$.  
\medskip\\  
For $f \in \frak{F}_{\alpha \beta}$ and $p \in Q(\theta _\alpha)$ we may define $\bar f[p]$ with $dom(\bar f[p])=\bar f[dom(p)]$ by setting
$$\bar f[p](\bar f(\eta))=\bar f_\eta \otimes {f}[p(\eta)]\quad \hbox{for all}\quad\eta \in dom(p)$$
where $\bar f$, $\bar f_\eta$ are as in lemma 2.6 and
$$\bar f_\eta \otimes {f}:\varphi_\eta \times \omega\eta \rightarrow \varphi_{\bar f(\eta)} \times {\omega\bar f(\eta)}, \ \langle \gamma,\omega\delta+n\rangle \mapsto \langle \bar f_\eta (\gamma),{\omega f}(\delta)+n\rangle$$ 
for all $n\in \omega$
$$\bar f_\eta \otimes {f}:(\varphi_\eta \times \omega\eta)\times 2 \rightarrow (\varphi_{\bar f(\eta)} \times {\omega\bar f(\eta)})\times 2,\ \langle x,\epsilon\rangle \mapsto \langle  \bar f_\eta \otimes {\bar f}(x),\epsilon \rangle.$$
In the same way we may define $\overline{\pi ^\prime _{st}}[p]$.
\bigskip\\
The reason why we use $\bar f$ instead of $f \in \frak{F}_{\alpha\beta}$ is that $f$ does not map the support of a condition correctly. For an example, consider the case $\beta=\alpha+1$ and let $f \in \frak{F}_{\alpha\beta}$ be right-branching. Let $\delta$ be the splitting point of $f$, i.e. $f(\delta)=\theta_\alpha$. Assume that $p\in Q(\theta_\alpha)$, $\delta \in dom(p)$ and $dom(p(\delta))\subseteq \varphi_\delta \times \omega\delta$. Let $f[p]$ be defined by $dom(f[p])=f[dom(p)]$ and $f[p](f(\eta)):=f_\eta \otimes f[p(\eta)]$ for all $\eta \in dom(p)$. We will show that $f[p] \notin Q(\theta_\beta)$. To do so, notice first that $f_\delta=f^\#(\delta) \circ \bar f_\delta$ by lemma 2.6 (1). However, $\bar f_\delta = id \upharpoonright \varphi_\delta$, because $f$ is right-branching with splitting point $\delta$. So $f_\delta=f^\#(\delta)$. Hence $f[p](\theta_\alpha)=f_\delta \otimes f[p(\delta)] = f^\#(\delta)[p(\delta)]$ because $dom(p(\delta))\subseteq \varphi_\delta \times \omega\delta$ and $f \upharpoonright \delta=id \upharpoonright \delta$. However, this contradicts the fact that all $q \in Q(\theta_\beta)$ are of the form $q=r^\ast \upharpoonright supp(r)$ for some $r \in P(\omega_3)$ because in this case $q(\theta_\alpha)\neq g[\bar q]$ for all $g \in \frak{G}_{\gamma\theta_\alpha}$, $\bar q \in P(\varphi_\gamma)$ and $\gamma < \theta_\alpha$ by the definition of the support of a condition. 
\smallskip\\
This problem does obviously not occur, if we consider $\bar f[p]$.
\bigskip\\
{\bf Lemma 4.3}
\smallskip\\
(a) If  $f \in \frak{F}_{\alpha \beta}$ and $p \in Q(\theta _\alpha)$, then $\bar f[p] \in Q(\theta_\beta)$.
\smallskip\\
(b) If $s \prec ^\prime t$ and $p \in Q(\nu(s)+1)$, then $\overline{\pi ^\prime_{st}}[p]\in Q(\nu(t)+1)$.
\smallskip\\
{\bf Proof:} 
Set $q:=\bar f[p]$. Let $dom(p) =\{ \alpha _1<\dots < \alpha _n\}$ and $dom(q)=\{ \beta _1 < \dots < \beta_n \} := \{ \bar f(\alpha _1) < \dots < \bar f(\alpha _n)\}$. By the definition of the support of a condition, all $\alpha _i$ are successor ordinals. And $f(\alpha _i-1)=\bar f(\alpha _i)-1$ by the definition of $\bar f$. Set $q(\beta_i-1)=e_{\beta_i-1}(q(\beta_i))$. Then it suffices to prove that there are functions $g_i \in \frak{G}_{\beta_i,\beta _{i+1} -1}$ such that 
\smallskip\\
(1) $q(\beta_{i+1}-1)=g_i[q(\beta_i)]$
\smallskip\\
(2) $q(\beta _i)\notin rng(\sigma_{\beta_i-1})$, $q(\beta _i)\notin P(\varphi_{\beta _i-1})$:
\smallskip\\
Since $p$ is a condition, there are functions $j_i \in \frak{G}_{\alpha _i,\alpha _{i+1}-1}$ such that 
$$p(\alpha _{i+1}-1)=j_i[p(\alpha _i)].$$ So we can set $$g_i=f_{\alpha _i,\alpha _{i+1}-1}(j_i) \circ f^\#(\alpha_i).$$ 
We need to check (1). We first prove that
$$f_{\alpha _{i+1}-1}\otimes { f}[e_{\alpha_{i+1}-1}(p(\alpha_{i+1}))]=e_{\beta_{i+1}-1}(q(\beta_{i+1})).$$
To see this, we use lemma 2.6 (4) which says
$$\forall \xi < \zeta \ \forall  b \in \frak{G}_{\xi \zeta}\ \ \bar f_\zeta \circ b =\bar f_{\xi\zeta}(b) \circ f_\xi.$$
Applying it for $\xi=\alpha_{i+1}-1$, $\zeta=\alpha_{i+1}$ and $b=id \upharpoonright \varphi_{\alpha_{i+1}-1}$, we get$$q(\beta_{i+1}) \upharpoonright (\varphi_{\beta_{i+1}-1} \times \omega(\beta_{i+1}-1))=\bar f_{\alpha_{i+1}}\otimes f[p(\alpha_{i+1})] \upharpoonright  (\varphi_{\beta_{i+1}-1} \times \omega(\beta_{i+1}-1))=$$
$$ =f_{\alpha_{i+1}-1 }\otimes f[p(\alpha_{i+1}) \upharpoonright (\varphi_{\alpha_{i+1}-1} \times \omega(\alpha_{i+1}-1))]$$
where the first equality holds by the definition of $q=\bar f[p]$.
\smallskip\\
Applying it for $\xi=\alpha_{i+1}-1$, $\zeta=\alpha_{i+1}$ and the splitting map $b$ of $\frak{G}_{\alpha_{i+1}-1,\alpha_{i+1}}$, we obtain
$$\bar f_{\xi\zeta}(b)^{-1}[q(\beta_{i+1}) \upharpoonright (\varphi_{\beta_{i+1}} \times \omega(\beta_{i+1}-1))]=$$
$$=\bar f_{\xi\zeta}(b)^{-1}[\bar f_{\alpha_{i+1}}\otimes f[p(\alpha_{i+1})] \upharpoonright (\varphi_{\beta_{i+1}} \times \omega(\beta_{i+1}-1))]=$$
$$(\bar f_{\xi\zeta}(b)^{-1}\circ \bar f_{\alpha_{i+1}}) \otimes  f[p(\alpha_{i+1}) \upharpoonright (\varphi_{\alpha_{i+1}} \times \omega(\alpha_{i+1}-1))]=$$
$$=(f_\xi \circ b^{-1}) \otimes  f[p(\alpha_{i+1}) \upharpoonright (\varphi_{\alpha_{i+1}} \times \omega(\alpha_{i+1}-1))]=$$
$$=(f_\xi \otimes f)[b^{-1}[p(\alpha_{i+1}) \upharpoonright (\varphi_{\alpha_{i+1}} \times \omega(\alpha_{i+1}-1))]].$$
However, by definition
\smallskip\\
$e_{\beta_{i+1}-1}(q(\beta_{i+1}))=$
$$ q(\beta_{i+1}) \upharpoonright (\varphi_{\beta_{i+1}-1} \times \omega(\beta_{i+1}-1))\quad \cup \quad 
\bar f_{\xi\zeta}(b)^{-1}[q(\beta_{i+1}) \upharpoonright (\varphi_{\beta_{i+1}} \times \omega(\beta_{i+1}-1))]$$
and 
\smallskip\\
$e_{\alpha_{i+1}-1}(p(\alpha_{i+1}))=$
$$ p(\alpha_{i+1}) \upharpoonright (\varphi_{\alpha_{i+1}-1} \times \omega(\alpha_{i+1}-1))\quad \cup \quad
b^{-1}[p(\alpha_{i+1}) \upharpoonright (\varphi_{\alpha_{i+1}} \times \omega(\alpha_{i+1}-1))].$$
This proves that
$$f_{\alpha _{i+1}-1}\otimes { f}[e_{\alpha_{i+1}-1}(p(\alpha_{i+1}))]=e_{\beta_{i+1}-1}(q(\beta_{i+1})).$$
Hence 
$$q(\beta_{i+1}-1)=e_{\beta_{i+1}-1}(q(\beta_{i+1}))= f_{\alpha _{i+1}-1}\otimes { f}[e_{\alpha_{i+1}-1}(p(\alpha_{i+1}))]=$$
$$=f_{\alpha _{i+1}-1}\otimes { f}[j_i(p(\alpha_{i}))]=(f_{\alpha _{i+1}-1} \circ j_i)\otimes { f}[p(\alpha_{i})]=$$ $$= 
(f_{\alpha _i,\alpha _{i+1}-1}(j_i) \circ f_{\alpha_i})\otimes {\bar f}[p(\alpha _i)]$$ by (6) in the definition of embeddings. However, $f_{\alpha _i}=f^\#(\alpha_i)\circ \bar f_{\alpha_i}$ by lemma 2.6. So 
$$(f_{\alpha _i,\alpha _{i+1}-1}(j_i) \circ f_{\alpha_i})\otimes { f}[p(\alpha _i)]=$$
$$=(f_{\alpha _i,\alpha _{i+1}-1}(j_i) \circ f^\#(\alpha_i)\circ \bar f_{\alpha_i} )\otimes { f}[p(\alpha _i)]=$$
$$=f_{\alpha _i,\alpha _{i+1}-1}(j_i) \circ f^\#(\alpha_i)[ \bar f_{\alpha_i} \otimes { f}[p(\alpha _i)]]=$$
$$=f_{\alpha _i,\alpha _{i+1}-1}(j_i) \circ f^\#(\alpha_i)[q(\beta_i)]$$ and we are done. 
\smallskip\\
To see (2), notice that by the definition of the support of a condition $p(\alpha _i) \notin rng(\sigma_{\alpha_i-1})$ and $p(\alpha_i) \notin P(\varphi_{\alpha_i-1})$. Now, we can use lemma 2.6 (4) to obtain that $q(\beta _i) \notin rng(\sigma_{\beta_i-1})$ and $q(\alpha_i) \notin P(\varphi_{\beta_i-1})$. The argument is very similar to the one we used to prove $$f_{\alpha _{i+1}-1}\otimes { f}[e_{\alpha_{i+1}-1}(p(\alpha_{i+1}))]=e_{\beta_{i+1}-1}(q(\beta_{i+1})).\quad \Box$$
\bigskip\\
In the following we thin out $Q(\gamma)$ to $\Q _\gamma$ to obtain a FS system along our gap-2 morass.
\smallskip\\
We define $\Q _\gamma$ by induction on the levels of $\langle \langle \theta _\alpha \mid \alpha \leq \omega_1 \rangle , \langle \frak{F}^\prime_{\alpha\beta} \mid \alpha < \beta \leq \omega_1 \rangle \rangle$.
\medskip\\
{\it Base Case}: $\beta =0$
\medskip\\
Then we only need to define $\Q _1$.
\smallskip\\
Let $\Q_1=Q(1)$.  
\medskip\\
{\it Successor Case}: $\beta = \alpha +1$
\medskip\\
We first define $\Q _{\theta _\beta}$. To do so, let $\P_{\varphi_{\theta_\beta}}$ be the set of all $p \in P(\varphi_{\theta_\beta})$ such that 
\smallskip\\
(1) $(\bar h_{\theta_\alpha}\otimes \bar h)^{-1}[p] \in \P_{\varphi_{\theta_\alpha}}$
\smallskip\\ 
(2) $g^{-1}[p\upharpoonright (\varphi_{\theta_\beta} \times \theta_\alpha)]$ and $(\bar h_{\theta_\alpha}\otimes \bar h)^{-1}[p]$ are compatible for all $g\in \frak{G}_{\theta_\alpha\theta_\beta}$
\smallskip\\
where $h$ is the unique right-branching embedding of $\frak{F}_{\alpha\beta}$.  
\bigskip\\
Set
$$\Q_{\theta_\beta}=\{ p^\ast \upharpoonright (supp(p)\cap\theta_\beta) \mid p \in \P_{\varphi_{\theta_\beta}}\} .$$
\smallskip\\
For $t \in T^\prime_\beta$ set $\Q_{\nu (t)+1}=\{ p \in \Q_{\theta_\beta} \mid dom(p) \subseteq \nu(t)+1\}$ and $\Q _\lambda =\bigcup \{ \Q_\eta \mid \eta < \lambda\}$ for $\lambda \in Lim$. 
\medskip\\
Set 
$$\sigma^\prime _{st}:\Q_{\nu (s)+1} \rightarrow \Q_{\nu(t)+1} ,\ p \mapsto \overline{\pi^\prime _{st}}[p].$$ 
It remains to define $e^\prime_\alpha$. If $p \in rng(\sigma^\prime _\alpha)$, then set $e^\prime_\alpha (p)=\sigma^{\prime-1}_\alpha(p)$. If $p \in \Q_{\theta _\alpha}$, then set $e^\prime_\alpha(p)=p$. And if $p \notin rng(\sigma^\prime_\alpha) \cup  \Q_{\theta _\alpha}$, then choose a $r \in \P_{\varphi_{\theta_\beta}}$ with $p=r^\ast \upharpoonright supp(r)$ and set
$$q:= \bigcup \{ g^{-1}[r \upharpoonright (\varphi_{\theta_\beta} \times \theta_\alpha)] \mid g\in \frak{G}_{\theta_\alpha\theta_\beta}\} \cup (\bar h_{\theta_\alpha}\otimes \bar h)^{-1}[r]$$
$$=r^\ast(\theta_\alpha) \cup (\bar h_{\theta_\alpha}\otimes \bar h)^{-1}[r].$$
Set $e^\prime_\alpha(p)=q^\ast \upharpoonright (supp(q) \cap \theta_\alpha)$. 
\medskip\\
{\it Limit Case}: $\beta \in Lim$
\medskip\\
For $t \in T^\prime_\beta$ set $\Q_{\nu(t)+1}=\bigcup \{ \sigma^\prime _{st}[\Q_{\nu (s)+1}] \mid s \prec^\prime t \}$ and $\Q_\lambda =\bigcup \{\Q_\eta \mid \eta < \lambda\}$ for $\lambda \in Lim$ where $\sigma^\prime _{st}:\Q_{\nu (s)+1} \rightarrow \Q_{\nu(t)+1}, \ p \mapsto \overline{\pi^\prime _{st}}[p]$.
\bigskip\\
Finally, set $\P_\eta=\{p \in P(\eta) \mid p^\ast \upharpoonright supp(p) \in \Q _{\omega _2}\}$ and $\P:=\P_{\omega_3}$.
\bigskip\\
We think that some explanations are appropriate. Let us first compare our definition to Velleman's construction in \cite{Velleman1987a}. His proof of the gap-3 theorem is theorem 5.3 of \cite{Velleman1987a}. He has to construct a structure $\frak{A}$.  Assume that his $\kappa^+=\omega_1$. Then he constructs $\frak{A}$ by constructing for every $\alpha < \omega_1$ a structure $\frak{A}_\alpha$ and taking a direct limit. However, the system of elementary embeddings he uses to take the direct limit is not a linear commutative system. That is, we do not have for every $\alpha < \omega_1$ a single elementary embedding   $f:\frak{A}_\alpha\rightarrow \frak{A}$ but an elementary embedding $f^\ast :\frak{A}_\alpha\rightarrow \frak{A} $ for every $f\in \frak{F}_{\alpha\omega_1}$. Moreover, he has to require that his structures $\frak{A}_\alpha$ "mirror" the structure of  $\langle \langle \varphi _\zeta \mid \zeta \leq \theta_\alpha , \langle \frak{G}_{\zeta\xi} \mid \zeta < \xi \leq \theta_\alpha\rangle\rangle$. Similarly, we obtain $\P$ as the direct limit of the $\P_{\varphi_{\theta_\alpha}}$, which is shown in the next lemma. Moreover, we proceed in such a way that $\P_{\varphi_{\theta_\alpha}} \subseteq P(\varphi_{\theta_\alpha})$. Hence also our $\P_{\varphi_{\theta_\alpha}}$ "mirror" the structure of  $\langle \langle \varphi _\zeta \mid \zeta \leq \theta_\alpha , \langle \frak{G}_{\zeta\xi} \mid \zeta < \xi \leq \theta_\alpha\rangle\rangle$.  As in the case of Velleman's construction, this is necessary to define $\P_{\varphi_{\alpha+1}}$ in the successor step (cf. lemma 5.2 of \cite{Velleman1987a}). Let us make some further remarks.
\bigskip\\
{\tt Remark 1:} 
\smallskip\\
We postpone the proof that this definies indeed an FS system along our gap-2 morass $\frak{M}$. However, we check the crucial condition (FS$_2$6) already here. To do so, let $p \in \Q_{\theta_\beta}$ and $\beta=\alpha+1$. Let $r \in \P_{\varphi_{\theta_\beta}}$ be such that $p=r^\ast \upharpoonright supp(r)$ and
$$q:= r^\ast(\theta_\alpha) \cup (\bar h_{\theta_\alpha}\otimes  h)^{-1}[r]$$
where $h$ is the right-branching embedding of $\frak{F}_{\alpha\beta}$.
We have to prove that $s:=q^\ast \upharpoonright (supp(p) \cap \theta_\alpha) \in \Q_{\theta_\alpha}$ is a reduction of $p$ with respect to $\sigma^\prime_\alpha$ and $id \upharpoonright \Q_{\theta_\alpha}$. To do so, let $t \in \Q_{\theta_\alpha}$ with $t \leq s$. We have to find an $u \in \Q_{\theta_\beta}$ such that $u \leq p,\sigma^\prime_\alpha(t),t$. Notice first that
$$s\leq r^\ast \upharpoonright (supp(r) \cap \theta_\alpha)$$
and
$$ s \leq (\bar h_{\theta_\alpha} \otimes h)^{-1}[r]^\ast \upharpoonright supp((\bar h_{\theta_\alpha}\otimes h)^{-1}[r]).$$
Hence
$$t\leq r^\ast \upharpoonright (supp(r) \cap \theta_\alpha)$$
and
$$ t \leq (\bar h_{\theta_\alpha} \otimes h)^{-1}[r]^\ast \upharpoonright supp((\bar h_{\theta_\alpha}\otimes h)^{-1}[r]).$$
Let $\nu =max(dom(t))$. Then $t(\nu)$ and $q$ are compatible. Set $v =q \cup t(\nu) \in \P_{\varphi_{\theta_\alpha}}$ and $w=r \cup v \cup (\bar h_{\theta_\alpha} \otimes h)[v]$. Then $w \leq r,t(\nu),(\bar h_{\theta_\alpha} \otimes h)[t(\nu)]$. Hence $u:= w^\ast \upharpoonright supp(w) \leq p$ because $w \leq r$. Moreover, $u \leq t,\sigma^\prime_\alpha(t)$. This is proved from $w \leq t(\nu),(\bar h_{\theta_\alpha} \otimes h)[t(\nu)]$ as in the proof of
$$p \in \Q_{\theta_\alpha} \quad \wedge \quad f \in \frak{F}_{\alpha\beta} \quad \Rightarrow \quad \bar f[p] \in \Q_{\theta_\beta}.$$
\bigskip\\
{\tt Remark 2:}
\smallskip\\
Suppose $p \in \P$ is given. Let $G$ be any generic filter with $p \in G$. Let $F= \bigcup \{ p \mid p \in G\}$. Then by (2) in the successor step of the construction, $F$ is not only already determined on $dom(p)$, but a lot more of $F$ is already determined. Set 
$$D=\{ n \in\omega \mid \exists \delta,\gamma \ \langle \gamma,\omega\delta+n\rangle \in dom(p)\}.$$
Then it will turn out that $F$ is at least not yet determined on 
$$\omega_3 \times \{ \omega\delta+n \mid n \in \omega - D,\delta \in \omega_2\}.$$ 
This will be used in lemma 4.6, which is the crucial step for proving that $\P$ adds a Hausdorff space.
\bigskip\\
{\tt Remark 3:}
\smallskip\\
Assume that $\beta=\alpha+1$ and that $h$ is the right-branching embedding of $\frak{F}_ {\alpha\beta}$. Let $p_1,p_2 \in \P_{\varphi_{\theta_\alpha}}$ be compatible and $g \in \frak{G}_{\theta_\alpha\theta_\beta}$. Then also $g[p_1]$ and $\bar h_{\theta_\alpha}\otimes h[p_2]$ are compatible, i.e. $g[p_1]$ and  $\bar h_{\theta_\alpha}\otimes h[p_2]$ agree on the common part of their domains. To prove this, let 
$$\langle \gamma,\eta \rangle \in dom(g[p_1])\cap dom( \bar h_{\theta_\alpha}\otimes h[p_2])$$
$$g(\langle \gamma_1,\eta_1\rangle)=\langle \gamma,\eta\rangle \qquad   \bar h_{\theta_\alpha}\otimes h (\langle \gamma_2,\eta_2\rangle)=\langle \gamma,\eta\rangle.$$
Since $h$ is right-branching, $\bar h_{\theta_\alpha}=h_{\theta_\alpha}$. Let $\delta$ be the critical point of $f \upharpoonright \theta_\alpha$. Then $\eta < \omega\delta$ and therefore $\eta=\eta_1=\eta_2$.
By (6) in the definition of right-branching, there exists a $b \in \frak{G}_{\delta\theta_\alpha}$ such that $f_{\delta\theta_\alpha}(b)=g$. Hence, by (6) in the definition of embedding,
$$h_{\theta_\alpha}\circ b=g \circ h_\delta.$$
So there exists $\langle \bar \gamma,\eta\rangle \in \varphi_\delta \times \omega\delta$ such that
$$h_{\theta_\alpha} \circ b(\langle \bar \gamma,\eta\rangle)=g \circ h_\delta(\langle \bar \gamma,\eta\rangle)=\langle \gamma,\eta\rangle$$
$$ h_\delta(\langle \bar \gamma,\eta\rangle)=\langle\gamma_1,\eta\rangle \qquad b(\langle \bar\gamma,\eta\rangle)=\langle \gamma_2,\eta\rangle.$$
By (5) in the definition of right-branching embedding, $h_\delta\in \frak{G}_{\delta\theta_\alpha}$. Hence $p_1(\gamma_1,\eta)=p^\ast_1(\delta)(\bar\gamma,\eta)$. Moreover, $p_2(\gamma_2,\eta)=p^\ast_2(\delta)(\bar\gamma,\eta)$ because $b \in \frak{G}_{\delta\theta_\alpha}$. However, $p_1$ and $p_2$ are compatible. Therefore, also $p^\ast_1(\delta)$ and $p^\ast_2(\delta)$ are compatible. So $p^\ast_1(\delta)(\bar\gamma,\eta)=p^\ast_2(\delta)(\bar\gamma,\eta)$. This in turn implies $p_1(\gamma_1,\eta)=p_2(\gamma_2,\eta)$. Hence $g[p_1](\gamma,\eta)=\bar h_{\theta_\alpha}\otimes h[p_2](\gamma,\eta)$. That's what we wanted to show.
\bigskip\\
The same argument shows for all $p \in\P_{\varphi_{\theta_\alpha}}$ and all $g \in \frak{G}_{\theta_\alpha\theta_\beta}$ that $g[p]\in \P_{\varphi_{\theta_\beta}}$, $\bar h_{\theta_\alpha}\otimes h[p] \in \P_{\varphi_{\theta_\beta}}$ and $g[p] \cup (\bar h_{\theta_\alpha}\otimes h)[p] \in \P_{\varphi_{\theta_\beta}}$.
\bigskip\\
For arbitrary $\alpha < \beta \leq \omega_1$ and $f \in \frak{F}_{\alpha\beta}$ define
$$f_{\theta_\alpha}\otimes f: \varphi_{\theta_\alpha} \times \omega\theta_\alpha \rightarrow \varphi_{\theta_\beta}\times \omega\theta_\beta, \quad \langle \gamma,\omega\delta +n\rangle \mapsto \langle f_{\theta_\alpha}(\gamma),\omega f(\delta)+n\rangle$$
for all $n \in \omega$ and
$$f_{\theta_\alpha}\otimes f: (\varphi_{\theta_\alpha} \times \omega\theta_\alpha)\times 2 \rightarrow (\varphi_{\theta_\beta}\times \omega\theta_\beta)\times 2, \quad \langle x,\epsilon\rangle \mapsto \langle f_{\theta_\alpha}\otimes f(x),\epsilon\rangle.$$
If $\beta=\alpha+1$, then $\frak{F}_{\alpha\beta}$ is an amalgamation by (3) in the definition of a simplified gap-2 morass. Hence $f \in \frak{F}_{\alpha\beta}$ is either left-branching or right-branching. Let $p \in \P_{\varphi_{\theta_\alpha}}$ and assume that $f$ is right-branching. Then $f_{\theta_\alpha}\otimes f[p]=\bar f_{\theta_\alpha}\otimes f[p]$ because $\bar f_{\theta_\alpha}=f_{\theta_\alpha}$. If $f$ is left-branching, then $f_{\theta_\alpha} \in \frak{G}_{\theta_\alpha\theta_\beta}$ and $f \upharpoonright \theta_\alpha=id \upharpoonright \theta_\alpha$. Hence $f_{\theta_\alpha}\otimes f[p]=f_{\theta_\alpha}[p]$. So in both cases
$$f_{\theta_\alpha}\otimes f[p]\in \P_{\varphi_{\theta_\beta}}.$$
By induction, this is also true if $\beta=\alpha+n$ for some $n \in \omega$. What does happen at limit levels?
\bigskip\\
{\bf Lemma 4.4}
\smallskip\\
For all  $\beta \in Lim$, $\P_{\theta_\beta}=\bigcup \{ f^\#(\theta_\alpha)\circ(\bar f_{\theta_\alpha}\otimes {\bar f})[\P_{\varphi_{\theta_\alpha}}] \mid  f \in \frak{F}_{\alpha\beta},\alpha < \beta\}.$
\smallskip\\
{\bf Proof:} We first prove $\supseteq$. Let $\alpha < \beta$, $p \in \P_{\varphi_{\theta_\alpha}}$ and $f \in \frak{F}_{\alpha\beta}$. We have to prove that $r:=f_{\theta_\alpha} \otimes f[p] \in \P_{\varphi_{\theta_\beta}}$. That is, we have to show that $r^\ast \upharpoonright supp(r) \in \Q_{\theta_\beta}$. But by the argument of lemma 4.3, $r^\ast \upharpoonright supp(r)=\bar f[q]$ where $q:= p^\ast \upharpoonright supp(p) \in \Q_{\theta_\alpha}$. Hence $\bar f[q]=r^\ast \upharpoonright supp(r)$ by the definition of $\Q_{\theta_\beta}$.
\smallskip\\
For the converse, let $p \in \P_{\varphi_{\theta_\beta}}$. Hence $r:=p \upharpoonright supp(p) \in \Q_{\theta_\beta}$ by the definition of $\P_{\varphi_{\theta_\beta}}$. 
Set $\nu := max(dom(r))$ and $t:=\langle \beta, \nu\rangle$. Moreover, let $g \in \frak{G}_{\nu\theta_\beta}$ be such that $g[r(\nu)]=p$. Let, by the definition of $\Q_{\nu+1}$, $s \prec^\prime t$ be such that $r=\sigma^\prime_{st}(\bar r)$ for some $\bar r$ in $\Q_{\nu(s)+1}$. Hence $r=\bar f[\bar r]$ for some $f \in \frak{F}_{\alpha\beta}$ such that $s:=\langle \alpha , \bar\nu\rangle$ and $f(\bar\nu)=\nu$. In particular, also $\bar f(\bar \nu)=\nu$. That is, if we set $\nu=\xi+1$, then $\xi \in rng(f)$. Hence $\bar f_{\bar\nu}=f_{\bar\nu}$ and $r(\nu)=f_{\bar\nu}\otimes f[r(\bar\nu)]$. Moreover, by (5)(c) in the definition of a simplified gap-2 morass, we may assume that $g=f_{\bar\nu\theta_\alpha}(\bar g)$ for some $\bar g \in \frak{G}_{\bar\nu\theta_\alpha}$. But then $p= f_{\theta_\alpha}\otimes { f}[\bar p]$ where $\bar p=\bar g [\bar r(\bar \nu)]$ by (6) in the definition of embedding. $\Box$
\bigskip\\
{\bf Lemma 4.5}
$$\langle\langle \P _\eta \mid \eta \leq \kappa ^{++} \rangle ,\langle \sigma _{st} \mid s \prec t \rangle ,\langle \sigma^\prime _{st} \mid s \prec^\prime t \rangle , \langle e_\alpha \mid \alpha < \kappa^+ \rangle, \langle e^\prime_\alpha \mid \alpha < \kappa \rangle\rangle$$ is a FS system along $\frak{M}$. Hence $\P_{\omega_3}$ is ccc.
\smallskip\\
{\bf Proof:} (FS$_2$2), (FS$_2$3), (FS$_2$4), (FS$_2$5) and (FS$_2$7) are clear from the construction. (FS$_2$6) was proved in remark 1. So we are only left with (FS$_2$1). That is, we have to prove that
$$\langle\langle \P _\eta \mid \eta \leq \kappa ^{++} \rangle ,\langle \sigma _{st} \mid s \prec t \rangle , \langle e_\alpha \mid \alpha < \kappa^+ \rangle\rangle$$ is a FS system along  $\langle \langle \varphi _\zeta \mid \zeta \leq \omega_2 \rangle , \langle \frak{G}_{\zeta\xi}\mid \zeta < \xi \leq \omega_2 \rangle\rangle$. We know that
 $$\langle\langle P (\eta) \mid \eta \leq \omega _3 \rangle ,\langle \sigma _{st} \mid s \prec t \rangle , \langle e_\alpha \mid \alpha < \omega_2 \rangle\rangle$$ is a FS system along  $\langle \langle \varphi _\zeta \mid \zeta \leq \omega_2 \rangle , \langle \frak{G}_{\zeta\xi}\mid \zeta < \xi \leq \omega_2 \rangle\rangle$. From this it follows immediately that (FS4), (FS5) and (FS7) also hold for
$$\langle\langle \P _\eta \mid \eta \leq \kappa ^{++} \rangle ,\langle \sigma _{st} \mid s \prec t \rangle , \langle e_\alpha \mid \alpha < \kappa^+ \rangle\rangle.$$
Moreover, (FS1) holds, because
$$(\ast) \qquad\P_\eta=\{ p \in \P \mid p \in P(\eta)\}$$
and for $P(\eta)$ we know (FS1) already. By $(\ast)$, one has to prove for (FS2), (FS3) and (FS6) that certain conditions are elements of $\P$. In the case of (FS2), for example, one has to show that $\sigma_{st}(p) \in \P_{\nu(t)+1}$ for all $p \in \P_{\nu(s)+1}$. In all three cases that's not difficult. $\Box$
\bigskip\\
The next two lemmas correspond to lemma 5.2 and lemma 5.3 of \cite{Irrgang4}. Lemma 4.6 will ensure that the generic topological space is Hausdorff. Lemma 4.7 will guarantee that the space has spread $\omega_1$.  
\bigskip\\
{\bf Lemma 4.6}
\smallskip\\
Let $p \in \P$ and $\gamma \neq \delta \in \omega_3$. Then there is $q \leq p$ in $\P$ and $\mu \in \omega_2$ such that $q(\gamma,\mu) \neq q(\delta,\mu)$.
\smallskip\\
{\bf Proof:} 
We prove by induction over the levels of the gap-2 morass, which we enumerate by $\beta \leq \omega_1$, the following 
\smallskip\\
{\it Claim:} Let $p \in \P_{\varphi_{\theta_\beta}}$ and $\gamma \neq \delta \in \varphi_{\theta_\beta}$. Then there is $q \leq p$ in $\P_{\varphi_{\theta_\beta}}$ and $\mu \in \omega\theta_\beta$ such that $q(\gamma,\mu)\neq q(\delta,\mu)$.
\bigskip\\
{\it Base Case:} $\beta=0$
\smallskip\\
Trivial.
\bigskip\\
{\it Successor Case:} $\beta=\alpha+1$
\bigskip\\
Let $h$ be the right-branching embedding of $\frak{F}_{\alpha\beta}$.
We consider four cases.
\bigskip\\
{\tt Case 1:} $\gamma,\delta \in rng(h_{\theta_\alpha})$
\smallskip\\
Let $p \in \P_{\varphi_{\theta_\beta}}$ be given, $h_{\theta_\alpha}(\bar \gamma)=\gamma$ and $h_{\theta_\alpha}(\bar \delta)=\delta$. Set $\bar p =(h_{\theta_\alpha}\otimes h)^{-1}[p] \cup p^\ast(\theta_\alpha)$. By the induction hypothesis, there exists a $\bar q \in \P_{\varphi_{\theta_\alpha}}$ and a $\bar \mu =\omega\bar\tau +n \in \omega\theta_\alpha$ ($n\in \omega$) such that $\bar q \leq \bar p$ and $\bar q (\bar \gamma,\bar \mu) \neq \bar q (\bar \delta,\bar \mu)$. Set 
$$q=p \cup (h_{\theta_\alpha}\otimes h)[\bar q]$$ 
and $\mu =\omega h(\bar \tau)+n$. Then $q \in \P_{\varphi_{\theta_\beta}}$ by remark 3, $q \leq p$ and $q(\gamma,\mu)=\bar q(\bar\gamma,\bar \mu)\neq \bar q(\bar\delta,\bar\mu)=q(\delta,\mu)$.
\bigskip\\
{\tt Case 2:} $\gamma,\delta \notin rng(h_{\theta_\alpha})$
\smallskip\\
We consider two subcases. Assume first that $\theta_\beta \notin Lim$. Then choose some $\mu \in [\omega(\theta_\beta-1),\omega\theta_\beta[$ such that $\mu \notin \{ \tau _2 \mid \exists \tau_1 \ \langle \tau_1,\tau_2\rangle \in dom(p)\}$. Set
$$q=p \cup \{ \langle \langle \gamma ,\mu \rangle ,0\rangle , \langle \langle \delta,\mu \rangle ,1\rangle \}.$$
By the choice of $\mu$, $q \in P(\varphi_{\theta_\beta})$. According to the case which we are in, $q^\ast (\theta_\alpha)=p^\ast(\theta_\alpha)$ and $(h_{\theta_\alpha}\otimes h)^{-1}[q]=(h_{\theta_\alpha}\otimes h)^{-1}[p]$. Hence $q^\ast$ and $(h_{\theta_\alpha}\otimes h)^{-1}[q]$ are compatible because $q^\ast$ and $(h_{\theta_\alpha}\otimes h)^{-1}[q]$ are compatible. So $q \in \P_{\varphi_{\theta_\beta}}$ and it is obviously as wanted. 
\smallskip\\
Now, suppose that $\theta_\beta \in Lim$. Assume w.l.o.g. that $\gamma < \delta$. Set $t=\langle \theta_\beta ,\delta\rangle$. Let $s \prec t$ be minimal such that $\gamma \in rng(\pi_{st})$. Let $s \in T_\eta$. Pick $\mu \in [\omega \cdot max(\theta_\alpha,\eta),\omega\theta_\beta[$ such that $\mu \notin \{ \tau_2 \mid \exists \tau_1 \ \langle \tau_1,\tau_2\rangle \in dom(p)\}$. Set
$$q=p \cup \{ \langle\langle \gamma,\mu\rangle ,0\rangle ,\langle\langle \delta,\mu\rangle,1\rangle\}.$$
As in the first subcase, $q \in P(\varphi_{\theta_\beta})$ by the choice of $\mu$. Also as in the first subcase, we can see that $q \in \P_{\varphi_{\theta_\beta}}$. Hence $q$ is as wanted.
\bigskip\\
{\tt Case 3:} $\gamma \in rng(h_{\theta_\alpha})$, $\delta \notin rng(h_{\theta_\alpha})$
\smallskip\\
Again, we consider two subcases. Assume first that $\theta_\beta \notin Lim$. Then choose $\mu \in [\omega (\theta_\beta -1),\omega\theta_\beta[$ such that $\mu \notin \{ \tau_2 \mid \exists \tau_1 \langle \tau_1,\tau_2 \rangle \in dom(p)\}$. Let $h_{\theta_\alpha}(\bar\gamma)=\gamma$ and $\tilde h(\bar\mu)=\mu$ where $\tilde h(\omega\tau+n)=\omega h(\tau)+n$.
\smallskip\\
Let 
$$\bar p=p^\ast(\theta_\alpha) \cup (h_{\theta_\alpha}\otimes h)^{-1}[p].$$
Then there exists by the previous lemma in $\P_{\varphi_{\theta_\alpha}}$ a $\bar q \leq \bar p$ such that $\langle \bar\gamma,\bar\mu\rangle \in dom(\bar q)$. Set
$$r=p \cup (h_{\theta_\alpha}\otimes h)[\bar q]$$
and
$$q=r \cup \{ \langle\langle \delta,\mu\rangle ,\epsilon\rangle\}$$
where $h_{\theta_\alpha}\otimes h[q](\gamma,\eta)\neq\epsilon \in 2$.
\smallskip\\
By the choice of $\mu$, $q \in P(\varphi_{\theta_\beta})$. By remark 3, $r \in \P_{\varphi_{\theta_\beta}}$. Hence $r^\ast(\theta_\alpha)$ and $(h_{\theta_\alpha}\otimes h)^{-1}[r]$ are compatible. According to the case which we are in, $q^\ast(\theta_\alpha)=r^\ast(\theta_\alpha)$ and $(h_{\theta_\alpha} \otimes h)^{-1}[q]=(h_{\theta_\alpha} \otimes h)^{-1}[r]$. So also $q \in \P_{\varphi_{\theta_\beta}}$. It is also as wanted.
\smallskip\\
Now, suppose that $\theta_\beta\in Lim$. Assume w.l.o.g. that $\gamma < \delta$. Set $t=\langle \theta_\beta,\delta\rangle$. Let $s \prec t$ be minimal such that $\gamma \in rng(\pi_{st})$. Let $s \in T_\eta$. Pick $\mu \in [\omega\cdot max(\theta_\alpha,\eta),\omega\theta_\beta[$ such that $\mu \notin\{ \tau_2 \mid \exists \tau_1 \ \langle \tau_1,\tau_2\rangle \in dom(p)\}$. Let $h_{\theta_\alpha}(\bar\gamma)=\gamma$ and $\tilde h(\bar\mu)=\mu$ where $\tilde h(\omega\tau+n)=\omega h(\tau)+n$ for all $n\in\omega$. Let
$$\bar p=p^\ast(\theta_\alpha) \cup (h_{\theta_\alpha}\otimes h)^{-1}[p].$$
From now on, proceed exactly as in the first subcase.
\bigskip\\
{\tt Case 4:} $\gamma \notin rng(h_{\theta_\alpha})$, $\delta \in rng(h_{\theta_\alpha})$
\smallskip\\
Like case 4.
\bigskip\\
{\it Limit Case:} $\beta \in Lim$
\smallskip\\
By a previous lemma, $\P_{\varphi_{\theta_\beta}}=\bigcup \{ (f_{\theta_\alpha}\otimes f)[\P_{\varphi_{\theta_\alpha}}] \mid \alpha < \beta , f\in \frak{F}_{\alpha\beta}\}.$ By (5) in the definition of a simplified gap-2 morass, $\varphi_{\theta_\beta}=\bigcup \{ f_{\theta_\alpha}[\varphi_{\theta_\alpha}] \mid \alpha < \beta, f \in \frak{F}_{\alpha\beta}\}$ and $\theta_\beta=\bigcup \{ f[\theta_\alpha]\mid \alpha < \beta ,f \in \frak{F}_{\alpha\beta}\}$. Hence by (4) in the definition of a simplified gap-2 morass, we can pick $\alpha < \beta$, $f \in \frak{F}_{\alpha\beta}$, $\bar p \in \P_{\varphi_{\theta_\alpha}}$, $\bar \gamma \in \varphi_{\theta_\alpha}$ and $\bar \delta \in \omega\theta_\alpha$ such that $f_{\theta_\alpha} \otimes f[\bar p]=p$, $f_{\theta_\alpha}(\bar\gamma)=\gamma$ and $\tilde f(\bar \delta)=\delta$ where $\tilde f(\omega\tau+n)=\omega f(\tau)+n$ for all $n \in \omega$. By the induction hypothesis, there exists $\bar q \leq \bar p$ such that $\bar q(\bar \gamma,\bar\mu)\neq\bar q(\bar \delta,\bar\mu)$. Set $q:= f_{\theta_\alpha}\otimes f[\bar q]$. Then $q$ is as wanted. $\Box$
\bigskip\\
{\bf Lemma 4.7}
\smallskip\\
Let $\langle p_i \mid i \in \omega_2\rangle$ be a sequence of conditions $p_i\in\P$ such that $p_i\neq p_j$ if $i\neq j$. Let $\langle\delta_i\mid i\in\omega_2\rangle$ be a sequence of ordinals $\delta_i\in\omega_3$ such that $\delta_i \in dom(x_{p_i})$ for all $i\in\omega_2$. Then there exist $i\neq j$ and $p\in\P$ such that $p\leq p_i,p_j$, $\langle\delta_i,\mu\rangle, \langle \delta_j,\mu\rangle \in x_p$ and $p(\delta_i,\mu)=p(\delta_j,\mu)$ for all $\mu \in rng(x_{p_j})$.
\smallskip\\
{\bf Proof:} By first extending the conditions, we may assume that $x_{p_i} = dom(x_{p_i}) \times rng(x_{p_i})$ for all $i \in \omega_2$. Hence $\langle\delta_j,\mu\rangle \in x_p$ will hold for all $\mu \in rng(x_{p_j})$ automatically. Moreover, we can assume by the $\Delta$-system lemma that all $x_{p_i}$ are isomorphic relative to the order of the ordinals, that $p_i \cong p_j$ for all $i,j \in \omega_2$, that $\pi(\delta_i)=\delta_j$ if $\pi:dom(x_{p_i})\cong dom(x_{p_j})$, that $\{ rng(x_{p_i}) \mid i \in \omega_2\}$ forms a $\Delta$-system with root $\Delta$, and that $\pi \upharpoonright \Delta=id\upharpoonright\Delta$ if $\pi:rng(x_{p_i}) \cong rng(x_{p_j})$. To prove the lemma, we consider two cases.
\bigskip\\
{\tt Case 1:} $rng(x_{p_i})=\Delta$ for all $i \in \omega_2$
\smallskip\\
Then we set $\eta=max(\Delta)$. Since there are $\omega_2$-many $p_i$ while $\P_{\varphi_{\eta+1}}$ has only $\omega_1$-many elements, there exist $p_i$ and $p_j$ with $i\neq j$ such that $p_i^\ast(\eta+1)=p_j^\ast(\eta+1)$. Hence by the usual arguments $p_i$ and $p_j$ are compatible. Set $p=p_i \cup p_j$. Then $p$ is as wanted, because $p_i \cong p_j$ and $\pi(\delta_i)=\delta_j$ if $\pi:dom(x_{p_i}) \cong dom(x_{p_j})$.
\bigskip\\
{\tt Case 2:} $rng(x_{p_i}) \neq \Delta$ for all $i \in \omega_2$
\smallskip\\
Then $\{ min(rng(x_{p_i}) -\Delta) \mid i \in \omega_2\}$ is unbounded in $\omega_2$. For every $i\in \omega_2$ choose $\alpha_i < \omega_1$, $f_i \in \frak{F}_{\alpha_i\omega_1}$, $\bar\delta_i \in \varphi_{\theta_{\alpha_i}}$ and $\bar p_i \in \P_{\varphi_{\theta_{\alpha_i}}}$ such that 
$$p_i=(f_i)_{\theta_{\alpha_i}}\otimes f_i[\bar p_i] \quad \hbox{and}\quad \delta_i=(f_i)_{\theta_{\alpha_i}}(\bar\delta_i).$$
Since there are $\omega_2$-many $\delta_i$ and $p_i$ but only $\omega_1$-many possible $\bar \delta_i$ and $\bar p_i$, we can assume that $\alpha_i=\alpha_j$, $\bar\delta_i=\bar\delta_j$ and $\bar p_i=\bar p_j$ for all $i,j \in \omega_2$. Set $\bar p=\bar p_i$, $\alpha=\alpha_i$ and $\bar \delta=\bar\delta_i$.
Let $\nu \in\omega_3$ be such that $p_i \in \P_\nu$ for all $i \in \omega_2$. Let $t=\langle \omega_2,\nu\rangle$. Let $s \prec t$ such that $p_i \in rng(\sigma_{st})$ for $\omega_1$-many $i \in \omega_2$.
Let $s\in T_\eta$. Pick $p_i$ such that $min(rng(x_{p_i})-\Delta)>\omega\eta$.
Let $\eta_i=min(rng(x_{p_i})-\Delta)$. Then by the choice of $f_i$, $\eta_i \in rng(f_i\upharpoonright \theta_\alpha)$. Let $u \prec t$ be such that $u \in T_{\eta_i}$. Let $f_i(\bar \eta_i)=\eta_i$.
Since there are $\omega_1$-many $j \in  \omega_2$ such that $p_j \in rng(\sigma_{st})$, there are also $\omega_1$-many $j \in \omega_2$ such that $p_j \in rng(\sigma_{ut})$. On the other hand, $rng((f_i)_{\bar\eta_i})$ is countable. So we can pick a $j \in \omega_2$ such that $\delta \notin rng((f_i)_{\bar\eta_i})$, $\pi_{ut}(\delta)=\delta_j$ and $p_j \in rng(\sigma_{ut})$. In the following we will show that there exists $p \leq p_i,p_j$ such that $\langle \delta_j,\mu\rangle \in x_p$ and $p(\delta_i,\mu)=p(\delta_j,\mu)$ for all $\mu \in rng(x_{p_i})$.
\medskip\\
For $\alpha < \beta \leq \omega_1$, let $f_i=g_i^\beta \circ j_i^\beta$ where $g_i^\beta \in \frak{F}_{\alpha\beta}$ and $j_i^\beta \in \frak{F}_{\beta\omega_1}$. Let $g_i^\beta(\eta_i^\beta)=\eta_i$ and $\gamma$ be minimal such that $\delta \in rng((g^\gamma_i)_{\eta_i^\gamma})$.
For $\gamma \leq \beta \leq \omega_1$, let $(g^\beta_i)_{\eta_i^\beta}(\delta^\beta)=\delta$, $p_i^\beta=(j_i^\beta)_{\theta_\alpha}\otimes j^\beta_i[\bar p]$, $g_i^\beta[\Delta_\beta]=\Delta$ and $\delta^\beta_i=(j^\beta_i)_{\theta_\alpha}(\bar\delta)$. We prove by induction over $\gamma \leq \beta \leq \omega_1$ the following
\bigskip\\
{\tt Claim 1:} If $\langle \eta^\beta_i,\delta^\beta\rangle \prec \langle\theta_\beta,\delta^\prime\rangle$, then there exists $p^\beta \leq p^\beta_i$ such that $\langle\delta^\prime,\mu\rangle\in x_{p^\beta}$ and $p^\beta(\delta_i^\beta,\mu)=p^\beta(\delta^\prime,\mu)$ for all $\mu \in rng(x_{p_i^\beta})-\Delta_\beta$.
\medskip\\
{\it Base case:} $\beta=\gamma$
\smallskip\\
By the definition of $\gamma$ and (5) in the definition of a simplified gap-2 morass, $\gamma$ is a successor ordinal. Let $\gamma=\gamma^\prime+1$. Moreover, $\theta_{\gamma^\prime} \leq \eta^\gamma_i$. Hence $p^\gamma_i=h_{\theta_{\gamma^\prime}}\otimes h[p_i^{\gamma^\prime}]$ where $h$ is the right-branching embedding of $\frak{F}_{\gamma^\prime\gamma}$.
We first notice, that $\delta^\prime \notin rng(h_{\theta_{\gamma^\prime}})$. Assume that this was not the case. Then pick a $\pi \in \frak{G}_{\eta^\gamma_i\theta_\gamma}$ such that $\pi(\delta^\beta)=\delta^\prime$. By (6) in the definition of right-branching, there is a $\bar\pi \in \frak{G}_{\eta_i^{\gamma^\prime}\theta_{\gamma^\prime}}$ such that $h_{\eta_i^{\gamma^\prime}\theta_{\gamma^\prime}}(\bar\pi)=\pi$. Let $h_{\theta_{\gamma^\prime}}(\bar\delta^\prime)=\delta^\prime$. Let $\langle\eta_i^{\gamma^\prime},\rho\rangle \prec \langle \theta_{\gamma^\prime},\bar\delta^\prime\rangle$.
By (6) in the definition of embedding,
$$h_{\eta^{\gamma^\prime}_i} \circ \bar\pi=\pi \circ h_{\theta^{\gamma^\prime}}.$$
Hence $h_{\eta^{\gamma^\prime}_i}(\rho)=\delta^\beta$, which contradicts the definition of $\gamma$.
We can define a condition $p^\beta \leq p^\beta_i$, $p^\beta \in \P$ by setting
$$p^\beta=p_i^\beta \cup \{ \langle \langle \delta^\prime,\mu\rangle ,p^\beta_i(\delta_i^\beta,\mu)\rangle \mid \mu \in rng(x_{p^\beta_i})-\theta_{\gamma^\prime}\}.$$
This $p^\beta$ is as wanted.
\bigskip\\
{\it Successor step:} $\beta=\rho+1$
\medskip\\
We consider two cases:
\medskip\\
{\tt Case 1:} $p_i^\beta=g[p^\rho_i]$ for some $g \in \frak{G}_{\theta_\rho\theta_\beta}$
\medskip\\
In this case $\eta_i^\rho=\eta_i^\beta < \theta_\rho$. Let $\langle \eta_i^\beta,\delta^\beta\rangle \prec \langle\theta_\rho,\delta^{\prime\prime}\rangle \prec \langle\theta_\beta,\delta^\prime\rangle$. Let $\pi \in \frak{G}_{\theta_\rho\theta_\beta}$ such that $\pi(\delta^{\prime\prime}) =\delta^\prime$. Then by the induction hypothesis, there exists $p^\prime \leq p^\rho_i$ such that
$$p^\prime(\delta^{\prime\prime},\mu)=p^\prime(\delta^\rho_i,\mu)$$ 
for all $\mu \in rng(x_{p^\rho_i})-\Delta_\rho$. Set 
$$p^\beta =\pi [p^\prime] \cup g[p^\prime].$$ Then by remark 3, $p^\beta \in \P$ and
$$p^\beta(\delta^\prime,\mu)=p^\prime(\delta^{\prime\prime},\mu)=p^\prime(\delta^\rho_i,\mu)=p^\beta_i(\delta_i^\beta,\mu)$$
for all $\mu \in rng(x_{p_i^\rho})-\Delta_\rho= rng(x_{p_i^\beta})-\Delta_\beta$. Hence $p^\beta$ is as wanted.
\bigskip\\
{\tt Case 2:} $p^\beta_i=h_{\theta_\rho}[p^\rho_i]$ where $h$ is the right-branching embedding of $\frak{F}_{\rho\beta}$
\medskip\\
We consider three subcases.
\medskip\\
{\tt Subcase 1:} $\delta^\prime \in rng(h_{\theta_\rho})$
\medskip\\
Let $h_{\theta_\rho}(\bar\delta^\prime)=\delta^\prime$. Then by (6) in the definition of embedding, $\langle\eta^\rho_i,\delta^\rho\rangle \prec \langle\theta_\rho,\bar\delta^\prime\rangle$. 
Hence by the induction hypothesis, there exists $p^\rho \leq p_i^\rho$ such that $\langle\bar\delta^\prime,\mu\rangle \in x_{p^\rho}$ and $p^\rho(\delta_i^\rho,\mu)=p^\rho(\bar\delta^\prime,\mu)$ for all $\mu \in rng(x_{p^\rho_i})-\Delta_\rho$. Set 
$$p^\beta=h_{\theta_\rho} \otimes h[p^\rho].$$
Then $p^\beta$ is as wanted.
\bigskip\\
{\tt Subcase 2:} $\delta^\prime \notin rng(h_{\theta_\rho})$ and $\theta_\rho \leq \eta^\beta_i$
\medskip\\
Exactly like the base case of the induction.
\bigskip\\
{\tt Subcase 3:} $\delta^\prime \notin rng(h_{\theta_\rho})$ and $\eta^\beta_i < \theta_\rho$.
\medskip\\
This case is a combination of the base case of the induction and of case 1. Let $\langle\eta_i^\beta,\delta^\beta\rangle \prec \langle\theta_\rho,\delta^{\prime\prime}\rangle \prec \langle\theta_\beta,\delta^\prime\rangle$. Let $\pi \in \frak{G}_{\theta_\rho\theta_\beta}$ such that $\pi(\delta^{\prime\prime})=\delta^\prime$. Then by the induction hypothesis, there exists $p^\rho \leq p_i^\rho$ such that $\langle\bar\delta^\prime,\mu\rangle \in x_{p^\rho}$ and $$p^\rho(\delta_i^\rho,\mu)=p^\rho(\delta^{\prime\prime},\mu)$$ for all $\mu \in rng(x_{p^\rho_i})-\Delta_\rho$.
Set 
$$p^\beta=\pi[p^\rho] \cup (h_{\theta_\rho}\otimes h)[p^\rho] \cup\{ \langle\langle \delta^\prime,\mu\rangle ,p^\beta_i(\delta_i^\beta,\mu) \rangle \mid \mu \in rng(x_{p^\beta_i})-\theta_\rho\}.$$
By remark 3, $p^\beta \in \P$. We claim that $p^\beta$ is as wanted. For $\mu \in rng(x_{p_i^\beta})-\theta_\rho$,
$$p^\beta(\delta_i^\beta,\mu)=p^\beta(\delta^\prime,\mu)$$
holds by definition. For $\mu \in rng(x_{p_i^\beta}) \cap \theta_\rho=rng(x_{p^\rho_i})\cap \theta_\rho$, we have
$$p^\beta(\delta^\prime ,\mu)=p^\rho(\delta^{\prime\prime},\mu)=p^\rho(\delta_i^\rho,\mu)=p_i^\beta(\delta_i^\beta,\mu).$$
This finishes the proof of the successor step.
\medskip\\
{\it Limit case:} $\beta \in Lim$
\medskip\\
By lemma 4.4 and by (4) and (5) in the definition of a simplified gap-2 morass, we can pick a $\rho < \beta$ and a $f \in \frak{F}_{\rho\beta}$ such that $\delta^\prime\in rng(f_{\theta_\rho})$ and $f_{\theta_\rho}\otimes f[p^\rho_i]=p_i^\beta$. Let $f_{\theta_\rho}(\bar\delta^\prime)=\delta^\prime$. Then by (6) in the definition of embedding, $\langle\eta^\rho_i,\delta^\rho\rangle \prec \langle\theta_\rho,\bar\delta^{\prime}\rangle$.
Hence we can pick by the induction hypothesis a $p^\rho \leq p_i^\rho$ such that $\langle\bar\delta^\prime,\mu\rangle \in x_{p^\rho}$ and 
$$p^\rho(\delta_i^\rho,\mu)=p^\rho(\bar\delta^\prime,\mu)$$ 
for all $\mu \in rng(x_{p^\rho_i})-\Delta_\rho$. Set
$$p^\beta=f_{\theta_\rho}\otimes f[p^\rho].$$
Then $p^\beta$ is obviously as wanted. This finishes the proof of claim 1.
\bigskip\\
Finally, we can prove by induction over $\alpha < \beta \leq\omega_1$
\bigskip\\
{\tt Claim 2:} For $\alpha \leq \beta < \gamma$, set $p^\beta:=p^\beta_i$. For $\gamma \leq \beta < \omega_1$, let $p^\beta$ be as in claim 1. Then there exists for all $\alpha \leq \beta < \omega_1$ a $p \in \P$ such that $p \leq p^\beta,p_j^\beta$.
\medskip\\
{\it Base case:} $\beta=\alpha$
\smallskip\\
Trivial.
\medskip\\
{\it Successor case:} $\beta=\rho+1$
\medskip\\
We consider four cases.
\bigskip\\
{\tt Case 1:} $p_i^\beta=f[p^\rho_i]$ and $p_j^\beta=g[p_j^\rho]$ for some $f,g \in \frak{G}_{\theta_\rho\theta_\beta}$
\medskip\\
By the induction hypothesis, there exists a $\bar p \leq p^\rho,p_j^\rho$. Set
$$p=f[\bar p] \cup g[\bar p] \cup p^\beta.$$
It is not difficult to see that $p \in \P$ in all the different cases which occur in the definition of $p^\beta$.
\bigskip\\
{\tt Case 2:} $p_i^\beta=h_{\theta_\rho}\otimes h[p^\rho_i]$ and $p_j^\beta=g[p_j^\rho]$ where $g \in \frak{G}_{\theta_\rho\theta_\beta}$ and $h$ is the right-branching embedding of $\frak{F}_{\rho\beta}$
\medskip\\
By the induction hypothesis, there exists a $\bar p \leq p^\rho,p_j^\rho$. Set
$$p=g[\bar p] \cup (h_{\theta_\rho}\otimes h)[\bar p] \cup p^\beta.$$
It is not difficult to see that $p \in \P$ in all the different cases which occur in the definition of $p^\beta$.
\bigskip\\
{\tt Case 3:} $p_j^\beta=h_{\theta_\rho}\otimes h[p^\rho_j]$ and $p_i^\beta=g[p_i^\rho]$ where $g \in \frak{G}_{\theta_\rho\theta_\beta}$ and $h$ is the right-branching embedding of $\frak{F}_{\rho\beta}$
\medskip\\
Like case 2.
\bigskip\\
{\tt Case 4:} $p_i^\beta=h_{\theta_\rho}\otimes h[p^\rho_i]$ and $p_j^\beta=h_{\theta_\rho}\otimes h[p_j^\rho]$ where  $h$ is the right-branching embedding of $\frak{F}_{\rho\beta}$
\medskip\\
By the induction hypothesis, there exists a $\bar p \leq p^\rho,p_j^\rho$. Set
$$p= (h_{\theta_\rho}\otimes h)[\bar p] \cup p^\beta.$$
It is not difficult to see that $p \in \P$ in all the different cases which occur in the definition of $p^\beta$.
\bigskip\\
{\it Limit case:} $\beta \in Lim$
\medskip\\
This is proved very similar to the limit step in claim 1.
\bigskip\\
This finishes claim 2 and proves the lemma, if we set $\beta=\omega_1$ and $\delta^\prime=\delta_j$. $\Box$
\pagebreak\\
{\bf Lemma 4.8}
\smallskip\\
(a) $i:\P_{\omega_3} \rightarrow \Q_{\omega_2}, p \mapsto p^\ast \upharpoonright supp(p)$ is a dense embedding.
\smallskip\\
(b) There is a ccc-forcing $\bar \P$ of size $\omega_1$ such that $\Q_{\omega_2}$ embedds densely into $\bar \P$.
\smallskip\\
{\bf Proof:} (a) We have $i[\P_{\omega_3}]=\Q_{\omega_2}$. So it is clear, that $i[\P_{\omega_3}]$ is dense in $\Q_{\omega_2}$. It remains to check (1) and (2) of the definition of embedding. It follows from lemma 4.2, that (1) holds. For (2) assume first that $p,p^\prime \in \P_{\omega_3}$ are compatible. So there is $r \leq p,p^\prime$ in $\P_{\omega_3}$. Hence $i(r)\leq i(p),i(p^\prime)$ by lemma 4.2. So $i(p),i(p^\prime) \in \Q_{\omega_2}$ are compatible. Conversely assume that $i(p),i(p^\prime) \in \Q_{\omega_2}$ are compatible. Then $p,p^\prime \in \P_{\omega_3}$ are compatible by lemma 3.1. 
\smallskip\\
(b) Note, that $\langle\langle \Q_\eta \mid \eta \leq \omega_3 \rangle , \langle \sigma^\prime_{st} \mid s \prec ^\prime t \rangle , \langle e_\alpha^\prime \mid \alpha < \omega_1 \rangle\rangle$ is an FS system along $\langle \langle \theta_\alpha \mid \alpha \leq \omega_1\rangle , \langle \frak{F}^\prime _{\alpha\beta} \mid \alpha < \beta \leq \omega_1\rangle \rangle$. Hence we can define $\bar \P$ from $\Q_{\omega_2}$ like we defined $\Q_{\omega_2}$ from $\P_{\omega_3}$. That $\Q_{\omega_2}$ embedds densely into $\bar \P$ is proved like before. $\Box$  
\bigskip\\
Before we prove the main theorem, let us recall the definition of the spread of a topological space. Let $(X,\tau)$ be a topological space with topology $\tau$. A subset $D \subseteq X$ is called discrete if for every $x \in D$ there exists an $U \in \tau$ such that $U \cap D =\{ x\}$. The spread $s(X)$ of $X$ is defined as $s(X)=\omega \cdot sup\{ card(D) \mid D$ is a discrete subset of $X\}$. 
\bigskip\\
{\bf Theorem 4.9}
\smallskip\\
If there is an $(\omega_1,2)$-morass, then there is a ccc-forcing $\bar \P$ of size $\omega_1$ that adds a $0$-dimensional $T_2$ topology on $\omega_3$ which has spread $\leq\omega_1$.
\smallskip\\
{\bf Proof:} By lemma 4.8, $\P_{\omega_3}$ embedds densely into $\bar \P$. Hence $\P_{\omega_3}$ and $\bar \P$ yield the same generic extensions. So it suffices to prove that $\P:=\P_{\omega_3}$ adds a $0$-dimensional $T_2$ topology on $\omega_3$ which has spread $\omega_1$. By lemma 4.5, $\P$ is ccc. Therefore, it preserves cardinals. Let $G$ be $P$-generic. We set $F=\bigcup \{ p \mid p \in G\}$. Then $F:\omega_3 \times \omega_2 \rightarrow 2$ by a simple density argument. Let $\tau$ be the topology on $\omega_3$ generated by the sets $A^i_\nu:=\{\alpha \in \omega_3 \mid F(\alpha,\nu)=i\}$. Thus a base for $\tau$ is formed by the sets $B_\varepsilon:= \bigcap \{ A^{\varepsilon(\nu)}_\nu \mid \nu \in dom(\varepsilon)\}$ where $\varepsilon: dom(\varepsilon) \rightarrow 2$ is finite and $ dom(\varepsilon)\subseteq \omega_2$. Hence $\tau$ is $0$-dimensional. We claim that $\tau$ is as wanted.
\smallskip\\
We first show that it is $T_2$. We have to prove that for $\gamma \neq \delta$ there is some $\mu \in \omega_2$ such that $F(\gamma,\mu) \neq F(\delta,\mu)$. This is clear by the genericity of $G$ and lemma 4.6.
\smallskip\\
It remains to prove that $\tau$ has spread $\leq\omega_1$. Assume not.  Let $\dot X$, $\dot h$ and $\dot B$ be names  and $p \in \P$  a condition such that 
\smallskip\\
$p\Vdash (\dot X \subseteq \omega_3$, $\dot h: \omega_2 \rightarrow \dot X$ is bijective, $\dot B:\omega_2 \rightarrow V$, $\forall i \in \omega_2$ $\dot B(i)$ is a basic open set, $\forall i \neq j \in \omega_2 \ \dot h(i) \in \dot B(i) \wedge \dot h(j) \notin \dot B(i)$). 
\smallskip\\
For every $i \in \omega_2$ let $p_i \leq p$ and $\delta_i$, $\varepsilon_i$ be such that $p_i \Vdash \dot h(\check i)=\check \delta_i \wedge \dot B(i)=B_{\check \varepsilon_i}$. By the previous lemma, there are $i \neq j$ and $r \in \P$ such that $r \leq p_i,p_j$, $\langle \delta_i,\mu\rangle,\langle \delta_j,\mu\rangle \in x_r$ and $r(\delta_i,\mu)=r(\delta_j,\mu)$ for all $\mu \in rng(x_{p_j})$. Hence $r \Vdash \dot h(j) = \check \delta_j \in \dot B(i)$ which contradicts the definition of $p$.   
 $\Box$  
\bigskip\\
By a theorem of Hajnal and Juhasz \cite{HajnalJuhasz}, $card(X) \leq 2^{2^{s(x)}}$ for every Hausdorff space $X$ where $s(X)$ is its spread. By theorem 2.3, we can assume that GCH holds in the ground model where we construct our forcing. Since the forcing satisfies ccc and has size $\omega _1$, it preserves GCH by the usual argument for Cohen forcing. So in the generic extension $card(X)=2^{2^{s(x)}}$ holds for the generic space $X$. Hence the theorem answers Juhasz' question \cite{Juhasz}, if the second $exp$ is necessary in the case that $s(X)=\omega_1$. Moreover, the theorem of Hajnal and Juhasz shows that we cannot expect to be able to construct from an $(\omega_1,3)$-morass a ccc forcing of size $\omega_1$ which adds a $T_2$ space of size $\omega_4$ and spread $\omega_1$. If this was possible, we could find such a forcing in $L$. However, by the usual argument used for Cohen forcing it preserves GCH which contradicts the theorem of Hajnal and Juhasz. For similar reasons it is not possible to construct as in \cite{Irrgang4} along an $(\omega_1,2)$-morass a ccc forcing that adds an $\omega_3$-Suslin tree. There it is easier to see what goes wrong. Namely it is not possible to prove the necessary versions of lemma 5.2 and lemma 5.3 of \cite{Irrgang4}. This is prevented by condition (FS$_2$6) of the definition of a FS system along a gap-2 morass, which requires that an amalgamation of conditions like in remark 1 above is possible.   
\smallskip\\
On the other hand, the observation that $i:\P_{\omega_3} \rightarrow \Q_{\omega_2}, p \mapsto p^\ast \upharpoonright supp(p)$ is a dense embedding, also applies to the forcing which we constructed in \cite{Irrgang4}. This yields
\bigskip\\
{\bf Theorem 4.10}
\smallskip\\
If there is an $(\omega_1,1)$-morass, then there is ccc-forcing of size $\omega_1$ that adds an $\omega_2$-Suslin tree. $\Box$
\section{Local FS systems along morasses} 
In this section, we explain how the ideas from the previous sections can be used to construct forcings that can destroy GCH. As an example we reprove a consistency statement of Koszmider's \cite{Koszmider}. The same method can be used to construct ccc forcings that add an $(\omega,\omega_2)$-superatomic Boolean algebra or a witness for $\omega_2 \not\rightarrow (\omega:2)^2_\omega$. 
\smallskip\\
In the previous section, we observed that every forcing obtained by a FS system along a simplified $(\omega_1,1)$-morass preserves GCH, if lemma 4.2 holds for it and every $\P_\eta$ with $\eta < \omega_1$ is countable. However, these are exactly the most natural properties of forcings constructed by FS systems. So  all ``natural examples'' of FS systems along morasses seem to preserve GCH. So we can for example not expect to add a family $\{X_\alpha\mid \alpha < \omega_2\}$ of uncountable subsets $X_\alpha \subseteq \omega_1$ such that $X_\alpha \cap X_\beta$ is finite for any two $\alpha\neq \beta \in \omega_2$ because the existence of such a family implies $2^\omega \geq \omega_2$ by a result of Baumgartner's \cite{Baumgartner}. 
\smallskip\\
How can we overcome this difficulty? Can we obtain by a FS system along a $(\kappa,1)$-morass a normal, linear FS iteration $\P_{\kappa^+}$? Note, that then we automatically add $\kappa^+$-many new reals.
\medskip\\
Assume that $\langle \P_\eta \mid \eta \leq \kappa^+\rangle$ is a normal, linear FS iteration given as a set of $\kappa^+$-sequences $p \in \P_{\kappa^+}$ such that $\P_\eta=\{ p\upharpoonright \eta \mid p \in \P_{\kappa^+}\}$ and $\P_{\eta+1}\cong\P_\eta \ast \dot Q_\eta$ (where $\dot Q_\eta$ is a $\P_\eta$-name such that $\P_\eta \Vdash(\dot Q_\eta$ is a forcing)). Then $p:\kappa^+ \rightarrow V\in \P_{\kappa^+}$ iff $\P_\eta \Vdash p(\eta) \in \dot Q_\eta$ for all $\eta \in \kappa^+$ and $supp(p):=\{ \eta \in \kappa^+ \mid \P_\eta \not\Vdash p(\eta)=1_{\dot Q_\eta}\}$ is finite.
\medskip\\
For finite $\Delta \subseteq \kappa^+$ and $p \in \P_{\kappa^+}$ define $p_\Delta \in \P_{\kappa^+}$ by setting

$p_\Delta(\eta)=p(\eta)$ if $\eta \in \Delta$

$p_\Delta(\eta)=\dot 1_{Q_\eta}$ if $\eta \notin\Delta$

where $\dot 1_{Q_\eta}$ is a $\P_\eta$-name such that $\P_\eta \Vdash \dot 1_{Q_\eta}=1_{\dot Q_\eta}$.
\medskip\\
For $A\subseteq \P_{\kappa^+}$ and finite $\Delta \subseteq \kappa^+$ define
$$A_\Delta =\{ p_\Delta \mid p \in A \}.$$
If $\mu \geq \omega_1$ is regular and $\P_\Delta$ satisfies the $\mu$-cc for all finite $\Delta \subseteq \kappa^+$, then $\P_{\kappa^+}$ also satisfies the $\mu$-cc, as follows by the standard $\Delta$-system argument. 
\medskip\\
The idea is now to ensure the $\mu$-cc of every $\P_\Delta$ by constructing it by a FS system along a morass. This motivates the following definition: We say that a FS iteration $\langle \P_\eta \mid \eta \leq \kappa^+\rangle$ like above is a local FS system along a (simplified) $(\kappa,1)$-morass $\frak{M}$ iff for every finite $\Delta \subseteq \kappa^+$ there is a FS system $\langle\langle \Q^\Delta _\eta \mid \eta \leq \kappa ^{+} \rangle ,\langle \sigma^\Delta _{st} \mid s \prec t \rangle , \langle e^\Delta_\alpha \mid \alpha < \kappa \rangle\rangle$ along $\frak{M}$ such that $\P_\Delta \subseteq _\bot \Q^\Delta_{\kappa^+}$.  
\bigskip\\
So far, all this is of course only theory. As a simple example let me consider the forcing to add a chain $\langle X_\alpha \mid \alpha < \omega_2\rangle$ such that $X_\alpha \subseteq \omega_1$, $X_\beta -X_\alpha$ is finite and $X_\alpha-X_\beta$ has size $\omega_1$ for all $\beta < \alpha < \omega_2$. The natural forcing to do this would be
$$P:= \{ p:a_p \times b_p \rightarrow 2 \mid a_p \times b_p \subseteq \omega_2 \times \omega_1 \hbox{ finite }\}$$
where we set $p \leq q$ iff $q \subseteq p$ and
$$\forall \alpha_1 < \alpha_2 \in a_q \ \forall \beta \in b_p - b_q \ p(\alpha_1,\beta) \leq p(\alpha_2,\beta).$$
Obviously, we will set $X_\alpha=\{ \beta \in \omega_1 \mid p(\alpha, \beta)=1$ for some $p \in G\}$ for a $P$-generic $G$.
\smallskip\\
It is easily seen that $\langle P_\eta \mid \eta \leq \kappa^+\rangle$ with $P_\eta=\{ p \in P \mid a_p \subseteq \eta\}$ can be written as FS iteration such that $P_\Delta=\{ p \in P \mid a_p \subseteq \Delta\}$. On the other hand, it is not simply a product. Unfortunately, it also does not satisfy ccc. To see this, consider for every $\beta < \omega_1$ the function $p_\beta:\{0,1\} \times \{ \beta\} \rightarrow 2$ where $p_\beta(0,\beta)=1$ and $p_\beta(1,\beta)=0$. Then $A=\{ p_\beta \mid \beta \in \omega_1\}$ is an antichain of size $\omega_1$. Therefore, we need to thin out the forcing in an appropriate way. To do this, let $\langle \langle \theta _\alpha \mid \alpha \leq \omega_1 \rangle , \langle \frak{F}_{\alpha\beta}\mid \alpha < \beta \leq \omega_1 \rangle\rangle$ be a simplified $(\omega_1 ,1)$-morass. We will define a system $\langle\langle \P _\eta \mid \eta \leq \omega_2 \rangle ,\langle \sigma _{st} \mid s \prec t \rangle \rangle$ which satisfies properties (FS1) - (FS5) in the definition of FS system along a gap-1 morass.
\medskip\\
Let $\pi :\bar \theta \rightarrow \theta$ be a order-preserving map. Then $\pi :\bar\theta \rightarrow \theta$ induces maps $\pi : \bar\theta \times \omega_1 \rightarrow \theta \times\omega_1$ and  $\pi :(\bar \theta  \times \omega_1)\times 2  \rightarrow (\theta  \times \omega_1 ) \times 2$ in the obvious way:
$$\pi : \bar\theta \times \omega_1 \rightarrow \theta \times \omega_1,\quad \langle \gamma ,\delta\rangle \mapsto \langle \pi (\gamma ),\delta \rangle$$
$$\pi :(\bar \theta  \times \omega_1) \times 2 \rightarrow (\theta  \times \omega_1 ) \times 2, \quad \langle x ,\epsilon\rangle \mapsto \langle \pi (x ),\epsilon \rangle.$$
Basically we will define our maps $\sigma$ by setting $\sigma(p)=\pi[p]$.
\medskip\\
We define  $\langle\langle \P _\eta \mid \eta \leq \omega_2 \rangle ,\langle \sigma _{st} \mid s \prec t \rangle \rangle$ by induction on the levels of  $\langle \langle \theta _\alpha \mid \alpha \leq \omega_1 \rangle , \langle \frak{F}_{\alpha\beta}\mid \alpha < \beta \leq \omega_1 \rangle\rangle$ which we enumerate by $\beta \leq \omega_2$. 
\medskip\\
{\it Base Case}: $\beta =0$
\medskip\\
Then we need only to define $\P_1$.
\smallskip\\
Let $\P_1:=\{ p \in P \mid a_p \times b_p \subseteq 1\times 1\}$.  
\medskip\\
{\it Successor Case}: $\beta = \alpha +1$
\medskip\\
We first define $\P_{\theta_\beta}$. Let it be the set of all $p \in P$ such that:
\smallskip\\
(1) $a_p \times b_p \subseteq \theta_\beta \times \beta$.
\smallskip\\
(2) $f^{-1}_\alpha[p] \upharpoonright (\theta_\alpha \times \alpha) \in \P_{\theta_\alpha}$, $p \upharpoonright (\theta_\alpha \times \alpha) \in \P_{\theta_\alpha}$
where $h_\alpha$ is as in (P3) in the definition of a simplified gap-1 morass. 
\smallskip\\
(3) If $\alpha \in b_p$, then $p(\gamma,\alpha) \leq p(\delta,\alpha)$ for all $\gamma < \delta \in a_p$,
i.e.
$$p \upharpoonright (\theta_\beta \times \{ \alpha \}) \hbox{ is monotone.}$$
For all $\nu \leq \theta_\alpha$, $\P_\nu$ is already defined. For $\theta_\alpha <\nu \leq \theta_\beta$ set 
$$\P_{\nu}=\{ p \in \P_{\theta_\beta} \mid a_p\times b_p \subseteq \nu \times \beta \}.$$ 
Set 
$$\sigma _{st}:\P_{\nu (s)+1} \rightarrow \P_{\nu(t)+1} , p \mapsto \pi _{st}[p].$$ \smallskip\\
{\it Limit Case}: $\beta \in Lim$
\medskip\\
For $t \in T_\beta$ set $\P_{\nu(t)+1}=\bigcup \{ \sigma _{st}[\P_{\nu (s)+1}] \mid s \prec t \}$ and $\P_\lambda =\bigcup \{ \P_\eta \mid \eta < \lambda\}$ for $\lambda \in Lim$ where $\sigma _{st}:\P _{\nu (s)+1} \rightarrow P_{\nu(t)+1}, p \mapsto \pi _{st}[p]$.
\bigskip\\
Set $\P:=\P_{\omega_2}$.
\bigskip\\
A ccc forcing that adds a  chain $\langle X_\alpha \mid \alpha < \omega_2\rangle$ such that $X_\alpha \subseteq \omega_1$, $X_\beta -X_\alpha$ is finite and $X_\alpha-X_\beta$ has size $\omega_1$ for all $\beta < \alpha < \omega_2$ was first defined by Koszmider \cite{Koszmider}. He used Todorcevic's \cite{StevoBook} $\rho$-functions for his definition. In \cite{Morgan1996}, Morgan shows that it is possible to directly read off a $\rho$-function from a simplified gap-1 morass. If we use this $\rho$-function to define Koszmider's forcing, then we get exactly the same forcing  as with our approach.  
\bigskip\\
{\bf Lemma 5.1}
\smallskip\\
For $p \in P$, $p \in \P$ iff for all $\alpha < \omega_1$ and all $f \in \frak{F}_{\alpha+1,\omega_1}$
$$f^{-1}[p] \upharpoonright (\theta_{\alpha+1}\times\{\alpha\} )\quad \hbox{is monotone.}$$
{\bf Proof:} We prove by induction on $\gamma \leq \omega_1$ the following
\smallskip\\
{\it Claim:} $p \in \P_{\theta_\gamma}$ iff $p\in P$, $a_p \subseteq \theta_\gamma$, $b_p\subseteq\gamma$ and for all $\alpha < \gamma$ and all $f \in \frak{F}_{\alpha+1,\gamma}$
$$f^{-1}[p] \upharpoonright (\theta_{\alpha+1}  \times \{ \alpha\}) \hbox{ is monotone}.$$
{\it Base case:} $\gamma=0$
\smallskip\\
Then there is nothing to prove.
\medskip\\
{\it Successor case:} $\gamma=\beta+1$
\smallskip\\
Assume first that $p \in \P_{\theta_\gamma}$. Then, by (2) in the successor step of the definition of $\P_{\omega_2}$, $f^{-1}[p], (id \upharpoonright\theta_\beta)^{-1}[p] \in \P_{\theta_\beta}$. Now assume $f \in \frak{F}_{\alpha+1,\gamma}$ and $\alpha < \beta$. Then $f=f_\beta \circ f^\prime$ or $f=f^\prime$ for some $f^\prime \in \frak{F}_{\alpha+1,\beta}$ by (P2) and (P3). 
So by the induction hypothesis 
$$f^{-1}[p] \upharpoonright (\theta_{\alpha+1}  \times \{ \alpha\}) \hbox{ is monotone}$$
for all $f \in \frak{F}_{\alpha+1,\gamma}$ and all $\alpha < \beta$. Moreover, if $\alpha =\beta$ then the identity is the only $f\in \frak{F}_{\alpha+1,\gamma}$. In this case 
$$f^{-1}[p] \upharpoonright (\theta_{\alpha+1}  \times \{ \alpha\}) \hbox{ is monotone}$$
by (3) in the successor case of the definition of $\P$.
\smallskip\\
Now suppose that
$$f^{-1}[p] \upharpoonright (\theta_{\alpha+1}  \times \{ \alpha\}) \hbox{ is monotone}$$
for all $\alpha < \gamma$ and all $f \in \frak{F}_{\alpha+1,\gamma}$. We have to prove that (2) and (3) in the successor step of the definition of $\P$ hold. (3) obviously holds by the assumption because the identity is the only function in $\frak{F}_{\gamma\gamma}=\frak{F}_{\beta+1,\gamma}$. For (2), it suffices by the induction hypothesis to show that 
$$f^{-1}[h_\beta^{-1}[p]] \upharpoonright (\theta_{\alpha+1}  \times \{ \alpha\}) \hbox{ is monotone}$$
and
$$f^{-1}[(id\upharpoonright \theta_\beta)^{-1}[p]] \upharpoonright (\theta_{\alpha+1}  \times \{ \alpha\}) \hbox{ is monotone}$$
for all $f \in \frak{F}_{\alpha+1,\beta}$. This, however, holds by (P2) and the assumption.
\medskip\\
{\it Limit case:} $\gamma \in Lim$
\smallskip\\
Assume first that $p \in \P_{\theta_\beta}$. Let $\alpha < \gamma$ and $f \in \frak{F}_{\alpha+1,\gamma}$. We have to prove that
$$f^{-1}[p] \upharpoonright (\theta_{\alpha+1}  \times \{ \alpha\}) \hbox{ is monotone}.$$
By the limit step of the definition of $\P$, there are $\beta < \gamma$, $g \in \frak{F}_{\beta\gamma}$ and $\bar p\in \P_{\theta_\beta}$ such that $p=g[\bar p]$. By (P4) there are $\alpha,\beta<\delta<\gamma$, $g^\prime \in \frak{F}_{\beta\delta}$, $f^\prime\in \frak{F}_{\alpha\delta}$ and $j \in \frak{F}_{\delta\gamma}$ such that $g=j\circ g^\prime$ and $f=j\circ f^\prime$. Let $p^\prime:=g^\prime[\bar p]$. Then, by the induction hypothesis
$$(f^\prime)^{-1}[p^\prime] \upharpoonright (\theta_{\alpha+1}  \times \{ \alpha\}) \hbox{ is monotone}.$$
However, $(f^\prime)^{-1}[p^\prime]=(f^\prime)^{-1}[j^{-1}[p]]=f^{-1}$ and we are done. 
\smallskip\\
Now assume that
$$f^{-1}[p] \upharpoonright (\theta_{\alpha+1}  \times \{ \alpha\}) \hbox{ is monotone}$$
for all $\alpha < \gamma$ and all $f \in \frak{F}_{\alpha+1,\gamma}$. We have to prove that $p\in \P_{\theta_\gamma}$, i.e. that there exists $t\in T_\gamma$ and $s \prec t$ such that $p=\pi_{st}[\bar p]$ for some $\bar p \in \P_{\nu(s)+1}$. To find such $t$, $s \prec t$ and $\bar p$, let $\nu < \theta_\gamma$ be such that $a_p\subseteq \nu$. Since $\nu=\{\pi_{st}[\nu(s)] \mid s \prec t\}$ and $p:a_p\times b_p\rightarrow 2$ is finite, there exist $s \prec t$ such that $a_p\times b_p\subseteq rng(\pi_{st})$. Let $p=\pi_{st}[\bar p]$.
We need to prove that $\bar p \in \P_{\theta_\beta}$ where $\beta=\alpha(s)$. By the induction hypothesis  it suffices to prove that
$$f^{-1}[\bar p] \upharpoonright (\theta_{\alpha+1}  \times \{ \alpha\}) \hbox{ is monotone}$$
for all $\alpha < \beta$ and all $f\in \frak{F}_{\alpha+1,\beta}$. So let $f \in \frak{F}_{\alpha+1,\beta}$ and $g\in \frak{F}_{\beta\gamma}$ such that $\pi_{st}=g\upharpoonright \nu(s)+1$. Then  
$$f^{-1}[\bar p] \upharpoonright (\theta_{\alpha+1}  \times \{ \alpha\}) 
=f^{-1}[g^{-1}[p]] \upharpoonright (\theta_{\alpha+1}  \times \{ \alpha\})  =$$
$$=(g \circ f)^{-1}[p] \upharpoonright (\theta_{\alpha+1}  \times \{ \alpha\}) $$
which is monotone by our assumption. $\Box$
\bigskip\\
Unlike in the case of $\omega_2$-Suslin trees which we discussed in \cite{Irrgang4}, we cannot make $\langle\langle \P _\eta \mid \eta \leq \omega_2 \rangle ,\langle \sigma _{st} \mid s \prec t \rangle \rangle$ into a FS system along  $\langle \langle \theta _\alpha \mid \alpha \leq \omega_1 \rangle , \langle \frak{F}_{\alpha\beta}\mid \alpha < \beta \leq \omega_1 \rangle\rangle$ by adding an appropriate $\langle e_\alpha \mid \alpha < \omega_1\rangle$.
\medskip\\
Instead, we want to define for all finite $\Delta \subseteq \omega_2$ FS systems $\langle\langle \Q^\Delta _\eta \mid \eta \leq \omega_2 \rangle ,\langle \sigma^\Delta _{st} \mid s \prec t \rangle , \langle e^\Delta_\alpha \mid \alpha < \omega_1 \rangle\rangle$ along  $\langle \langle \theta _\alpha \mid \alpha \leq \omega_1 \rangle , \langle \frak{F}_{\alpha\beta}\mid \alpha < \beta \leq \omega_1 \rangle\rangle$ such that $\P_\Delta:=\{ p \in \P_{\omega_2} \mid a_p \subseteq \Delta\} \subseteq _\bot \Q^\Delta_{\omega_2}$. In other words, we want to represent every $p \in \P_\Delta$ as a function $p^\ast:\omega_1 \rightarrow V$ as in section 3 such that:
\smallskip\\
(1) $p^\ast(\alpha) \in \Q^\Delta_{\theta_\alpha}$ for all $\alpha < \omega_1$.
\smallskip\\
(2) If $p,q \in \P_\Delta$ and $p^\ast(\alpha),q^\ast(\alpha)$ are compatible in $\Q^\Delta_{\theta_\alpha}$ for $\alpha=max(supp(p) \cap supp(q))$, then $p$ and $q$ are compatible in $\P_\Delta$.
\medskip\\
How can we do this? Fix a finite $\Delta \subseteq \omega_2$. Set

$\eta=max(\Delta)$

$t=\langle \omega_1,\eta\rangle$

$s_0=min\{ s \prec t \mid \Delta \subseteq rng(\pi_{st})\}$

$\alpha_0 =\alpha(s_0)$.

Now, let $p \in \P_\Delta$. We simply set

$p^\ast(\alpha)=\pi^{-1}_{st}[p\upharpoonright (\omega_2 \times \alpha)]$ for $\alpha _0\leq\alpha < \omega_1$

where $s \in T_\alpha$, $s \prec t$. Like before we define

$supp(p)=\{ \alpha +1 \mid \alpha _0 \leq \alpha < \omega_1, p^\ast(\alpha +1)\neq p^\ast (\alpha), p^\ast(\alpha+1)\neq h_\alpha[p^\ast(\alpha)]\} \cup \{ \alpha_0\}$

where $h_\alpha$ is as in (P3) of the definition of a simplified gap-1 morass.
\smallskip\\
It is not completely obvious but easy to guess from this definition what the FS system  $\langle\langle \Q^\Delta _\eta \mid \eta \leq \omega_2 \rangle ,\langle \sigma^\Delta _{st} \mid s \prec t \rangle , \langle e^\Delta_\alpha \mid \alpha < \omega_1 \rangle\rangle$ looks like in the part above level $\alpha_0$. We could now explicitly give the definition of  $\langle\langle \Q^\Delta _\eta \mid \eta \leq \omega_2 \rangle ,\langle \sigma^\Delta _{st} \mid s \prec t \rangle , \langle e^\Delta_\alpha \mid \alpha < \omega_1 \rangle\rangle$ and infer from it that $\P_\Delta$ satisfies ccc. But this is very technical. Instead we will directly show the following
\bigskip\\
{\bf Lemma 5.2}
\smallskip\\
If $p,q \in \P_\Delta$ and $p^\ast(\alpha),q^\ast(\alpha)$ are compatible in $\P_{\theta_\alpha}$ for $\alpha=max(supp(p) \cap supp(q))$, then $p$ and $q$ are compatible in $\P_\Delta$.
\smallskip\\
{\bf Proof:} The proof is a simplified version of the proof of lemma 3.1. Suppose $p$ and $q$ are as in the lemma, but incompatible. Let $(supp(p)\cup supp(q))-\alpha =\{ \gamma_n < \dots < \gamma_1\}$. We prove by induction on $1 \leq i \leq n$, that $p^\ast(\gamma_i)$ and $q^\ast(\gamma_i)$ are incompatible for all $1 \leq i \leq n$. Since $\gamma_n=\alpha$, this yields the desired contradiction.
\smallskip\\
Note first, that $p^\ast(\gamma_1)$ and $q^\ast(\gamma_1)$ are incompatible because otherwise $p=\pi_{st}[p^\ast(\gamma_1)]$ and $q=\pi_{st}[q^\ast(\gamma_1)]$ were incompatible (for $s \in T_{\gamma_1}$, $s \prec t$). If $\gamma_1=\alpha$, we are done. So assume that $\gamma_1 \neq \alpha$. Then either $p^\ast(\gamma_1)=\pi_{\bar ss}[p^\ast(\gamma_1-1)]$ or $q^\ast(\gamma_1)=\pi_{\bar ss}[q^\ast(\gamma_1-1)]$ where $\bar s \prec s \prec t$, $\bar s \in T_{\gamma_1-1}$ and $s \in T_{\gamma_1}$. We assume in the following that $p^\ast(\gamma_1)=\pi_{\bar ss}[p^\ast(\gamma_1-1)]$. Mutatis mutandis, the other case works the same.

{\it Claim:} \quad $p^\ast(\gamma_1-1)$ and $q^\ast(\gamma_1-1)$ are incompatible in $\P_{\theta_{\gamma_1-1}}$

Assume not. Then there is $\bar r \leq p^\ast(\gamma_1-1),q^\ast(\gamma_1-1)$ in  $\P_{\varphi_{\gamma_1-1}}$ such that $a_{\bar r}=a_{p^\ast(\gamma_1-1)} \cup a_{q^\ast(\gamma_1-1)}$. Let $r^\prime:=\pi_{\bar ss}[\bar r]$. Then $r^\prime \leq \pi[p^\ast(\gamma_1-1)]=p^\ast(\gamma_1)$ and $r^\prime \leq \pi[q^\ast(\gamma_1-1)]=q^\ast(\gamma_1) \upharpoonright (\theta_{\gamma_1} \times \gamma_1)$. In the following we will construct an $r \leq p^\ast(\gamma_1),q^\ast(\gamma_1)$ which yields the contradiction we were looking for. By (2) in the definition of $\P_{\theta_{\gamma_1}}$, $\bar q(\eta,\gamma_1) \leq \bar q(\delta,\gamma_1)$ for all $\eta < \delta \in a_{\bar q}$ where $\bar q:= q^\ast(\gamma_1)$. Let $\tilde \delta = max \{ \delta \in a_{\bar q} \mid \bar q(\delta,\gamma_1)=0\}$ if the set is not empty. Otherwise, set $\tilde \delta=0$. Set
$$r=r^\prime \cup\{ \langle\langle \delta , \gamma_1\rangle ,0\rangle \mid \delta \leq \tilde \delta, \delta \in a_{r^\prime}\} \cup\{ \langle\langle \delta , \gamma_1\rangle ,1\rangle \mid \tilde \delta < \delta, \delta \in a_{r^\prime}\} .$$  
Then $r$ is as wanted. This proves the claim.

It follows from the claim, that $p^\ast(\gamma_2)$ and $q^\ast(\gamma_2)$ are incompatible. Hence we can prove the lemma by repeating this argument inductively finitely many times. $\Box$
\bigskip\\
{\bf Lemma 5.3}
\smallskip\\
$\P:=\P_{\omega_2}$ satisfies ccc.
\smallskip\\
{\bf Proof:} 
Let $A\subseteq \P$ be a set of size $\omega_1$. By the $\Delta$-lemma, we may assume that $\{ b_p \mid p \in A\}$ forms a $\Delta$-system with root $D$. We may moreover assume that for all $\alpha\in D$, all $f \in \frak{F}_{\alpha+1,\omega_1}$ and all $p,q \in A$
$$f^{-1}[p] \upharpoonright (\theta_{\alpha+1}  \times \{ \alpha\})  \subseteq f^{-1}[q] \upharpoonright (\theta_{\alpha+1}  \times \{ \alpha\})$$
or 
$$f^{-1}[p] \upharpoonright (\theta_{\alpha+1}  \times \{ \alpha\})  \supseteq f^{-1}[q] \upharpoonright (\theta_{\alpha+1}  \times \{ \alpha\}).$$
To see this assume that  $X=\{ a_p \mid p \in A\} \subseteq \omega_2$ forms a $\Delta$-system with root $\Delta_1$. Fix $\alpha \in D$. By thinning out $A$, we can ensure that whenever $a \neq b \in X$, $\eta \in a-b$, $\nu \in b-a$, $\eta < \nu$, $t= \langle\omega_1,\nu\rangle$, $s \prec t$, $s \in T_{\alpha+1}$, then $\eta \notin rng(\pi_{st})$. This suffices.
\smallskip\\
By the $\Delta$-system lemma, we may assume that $\{ a_p \mid p \in A\} \subseteq \omega_2$ forms a $\Delta$-system with root $\Delta_1$. Consider $A^\prime:=\{p \upharpoonright (\Delta_1 \times \omega_1)\mid p \in A\}$. By the $\Delta$-system lemma we may also assume that $\{ supp(p) \mid p \in A^\prime\} \subseteq \omega_1$ forms a $\Delta$-system with root $\Delta_2$. Let $\alpha=max(\Delta_2)$. Since $\P_{\theta_\alpha}$ is countable, there are $q_1 \neq q_2\in A^\prime$ such that $q^\ast_1(\alpha)=q_2^\ast(\alpha)$. Hence $q_1\neq q_2 \in A^\prime$ are compatible by a previous lemma. Assume that $q_1=p^\ast_1 \upharpoonright (\Delta_1\times\omega_1)$ and $q_2 = p^\ast_2 \upharpoonright (\Delta_1\times\omega_1)$ with $p_1,p_2 \in A$. 
We can define $p \leq p_1,p_2$ as follows: $a_p=a_{p_1}\cup a_{p_2}$, $b_p=b_{p_1}\cup b_{p_2}$, $p \upharpoonright (a_{p_1}\times b_{p_1})=p_1$,   $p \upharpoonright (a_{p_2}\times b_{p_2})=p_2$. We still need to define $p$ on $(a_p\times  b_p) -((a_{p_1} \times b_{p_1}) \cup (a_{p_2}\times b_{p_2}))$. We do this as in the previous lemma. That is, for $\beta \in b_p$ we set $\delta_\beta=max\{ \delta\in \Delta_1 \mid p(\delta,\beta)=0\}$ if this set is not empty. Otherwise, we set $\delta_\beta=0$. We set $p(\gamma,\beta)=1$ if we still need to define $p(\gamma,\beta)$ and $\gamma >\delta_\beta$. And we set $p(\gamma,\beta)=0$ if we still need to define $p(\gamma,\beta)$ and $\gamma \leq \delta_\beta$. Then $p \leq p_1,p_2$. We prove $p\leq p_1$. The other statement is showed similarly. Let $\gamma<\beta \in a_{p_1}$ and $\xi\in b_p-b_{p_1}$. We have to show that $p(\gamma,\xi)\leq p(\beta,\xi)$. If $\gamma,\beta\in \Delta_1$, then it holds because $q_1=p\upharpoonright(\Delta_1\times b_{q_1})$ and $q_2=p\upharpoonright(\Delta_1\times b_{q_2})$ are compatible. Otherwise, it holds by our definition of $p$ on  $(a_p\times  b_p) -((a_{p_1} \times b_{p_1}) \cup (a_{p_2}\times b_{p_2}))$.
\smallskip\\
It remains to prove that $p\in \P$. For this, we show that for all $\alpha <\omega_1$ and all $f\in \frak{F}_{\alpha+1,\omega_1}$
$$f^{-1}[p] \upharpoonright (\theta_{\alpha+1}  \times \{ \alpha\}) \hbox{ is monotone},$$
$$\hbox{i.e. } p \upharpoonright f[(\theta_{\alpha+1}  \times \{ \alpha\})]\hbox{ is monotone}.$$
Assume that $\alpha \in D$. Then by our second thinning-out
$$p_1 \upharpoonright f[(\theta_{\alpha+1}  \times \{ \alpha\})]  \subseteq p_2 \upharpoonright f[(\theta_{\alpha+1}  \times \{ \alpha\})]$$
or 
$$p_1 \upharpoonright f[(\theta_{\alpha+1}  \times \{ \alpha\})]  \supseteq p_2 \upharpoonright f[(\theta_{\alpha+1}  \times \{ \alpha\})]$$
and hence 
$$p \upharpoonright f[(\theta_{\alpha+1}  \times \{ \alpha\})]  = p_1 \upharpoonright f[(\theta_{\alpha+1}  \times \{ \alpha\})]\hbox{ is monotone}$$
or
$$p \upharpoonright f[(\theta_{\alpha+1}  \times \{ \alpha\})]  = p_2 \upharpoonright f[(\theta_{\alpha+1}  \times \{ \alpha\})]\hbox{ is monotone}.$$
Now, assume that $\alpha \notin D$. Then by our first thinning-out 
$$p_1 \upharpoonright f[(\theta_{\alpha+1}  \times \{ \alpha\}) ] =\emptyset \quad \hbox{or} \quad p_2 \upharpoonright f[ (\theta_{\alpha+1}  \times \{ \alpha\}) ] =\emptyset.$$
Hence
$$((a_p-a_{p_1})\times \{\alpha\})\cap f[ (\theta_{\alpha+1}  \times \{ \alpha\}) ]  = $$
$$=((a_p -(a_{p_1}  \cup a_{p_2}))\times \{\alpha\})\cap f[ (\theta_{\alpha+1}  \times \{ \alpha\}) ] $$
or
$$((a_p-a_{p_2})\times \{\alpha\})\cap f[ (\theta_{\alpha+1}  \times \{ \alpha\}) ]  = $$
$$=((a_p -(a_{p_1}  \cup a_{p_2}))\times \{\alpha\})\cap f[ (\theta_{\alpha+1}  \times \{ \alpha\}) ] .$$
To prove that $ p \upharpoonright f[(\theta_{\alpha+1}  \times \{ \alpha\})]\hbox{ is monotone},$ we consider the first case first. Let $\gamma<\delta \in f[\theta_{\alpha+1}]$. If $\gamma,\delta\in a_{p_1}$, then $p(\gamma,\alpha)=p_1(\gamma,\alpha)\leq p_1(\delta,\alpha)= p(\delta,\alpha)$ because $p_1 \in \P$. Otherwise $p(\gamma,\alpha)\leq p(\delta,\alpha)$ by the definition of $p$. The second case is proved in the same way where $p_1$ is replaced by $p_2$. $\Box$
\bigskip\\
{\bf Theorem 5.4}
\smallskip\\
If there is a simplified $(\omega_1,1)$-morass, then there is a ccc-forcing $\P$ which adds a chain $\langle X_\alpha \mid \alpha < \omega_2\rangle$ such that $X_\alpha \subseteq \omega_1$, $X_\beta -X_\alpha$ is finite and $X_\alpha-X_\beta$ has size $\omega_1$ for all $\beta < \alpha < \omega_2$. 
\smallskip\\
{\bf Proof:} By lemma 5.3, $\P$ satisfies ccc. Hence it preserves cardinals. It is easily seen by induction along the morass, that for every $\alpha \in \omega_2$ and every $\beta \in \omega_1$ the sets $D_\alpha=\{ p \in \P \mid \alpha \in a_p \}$ and $D^\prime_\beta=\{ p \in \P \mid \beta \in b_p\}$ are dense in $\P$. So if $G$ is $\P$-generic, then $F=\bigcup\{ p \mid p\in G\}$ is a function $F:\omega_2\times \omega_1 \rightarrow 2$. Set $X_\alpha =\{ \beta \in \omega_1 \mid F(\alpha,\beta)=1\}$. By the definition of $\leq $ on $\P$, $X_\beta - X_\alpha$ is finite for all $\beta < \alpha$. Finally, again by an easy induction along the morass we can prove that for all $\eta \in \omega_1$, $\beta < \alpha \in \omega_2$ the set $D^{\prime\prime}_{\eta,\alpha , \beta}=\{ p \in \P \mid \exists \gamma \geq \eta \ p(\beta , \gamma)=0 , p(\alpha , \gamma)=1\}$ is dense in $\P$. This yields that $X_\alpha -X_\beta$ is uncountable for all $\beta < \alpha < \omega_2$.  $\Box$

\bibliographystyle{plain}
\bibliography{biblio}

\begin{thebibliography}{10}

\bibitem{Baumgartner}
James~E. Baumgartner.
\newblock Almost-disjoint sets, the dense set problem and the partition
  calculus.
\newblock {\em Ann. Math. Logic}, 9(4):401--439, 1976.

\bibitem{BaumgartnerShelah}
James~E. Baumgartner and Saharon Shelah.
\newblock Remarks on superatomic {B}oolean algebras.
\newblock {\em Ann. Pure Appl. Logic}, 33(2):109--129, 1987.

\bibitem{Devlin}
Keith~J. Devlin.
\newblock {\em Constructibility}.
\newblock Perspectives in Mathematical Logic. Springer-Verlag, Berlin, 1984.

\bibitem{Donder}
Hans-Dieter Donder.
\newblock Another look at gap-{$1$} morasses.
\newblock In {\em Recursion theory (Ithaca, N.Y., 1982)}, volume~42 of {\em
  Proc. Sympos. Pure Math.}, pages 223--236. Amer. Math. Soc., Providence, RI,
  1985.

\bibitem{Fedorcuk}
V.~V. Fedor{\v{c}}uk.
\newblock The cardinality of hereditarily separable bicompacta.
\newblock {\em Dokl. Akad. Nauk SSSR}, 222(2):302--305, 1975.

\bibitem{Friedman}
Sy~D. Friedman.
\newblock {\em Fine structure and class forcing}, volume~3 of {\em de Gruyter
  Series in Logic and its Applications}.
\newblock Walter de Gruyter \& Co., Berlin, 2000.

\bibitem{HajnalJuhasz}
A.~Hajnal and I.~Juh{\'a}sz.
\newblock Discrete subspaces of topological spaces.
\newblock {\em Nederl. Akad. Wetensch. Proc. Ser. A 70=Indag. Math.},
  29:343--356, 1967.

\bibitem{Irrgang2}
Bernhard Irrgang.
\newblock Constructing $(\omega_1, \beta)$-morasses for $\omega _1 \leq \beta$.
\newblock Unpublished.

\bibitem{Irrgang1}
Bernhard Irrgang.
\newblock Proposing $(\omega_1, \beta)$-morasses for $\omega _1 \leq \beta$.
\newblock Unpublished.

\bibitem{Irrgang3}
Bernhard Irrgang.
\newblock Kondensation und {M}oraste.
\newblock Dissertation, Universit\"at M\"unchen, 2002.

\bibitem{Irrgang4}
Bernhard Irrgang.
\newblock Morasses and finite support iterations.
\newblock {\em Proc. Amer. Math. Soc.}, 137(3):1103--1113, 2009.

\bibitem{Jensen1}
Ronald~B. Jensen.
\newblock Higher-gap morasses.
\newblock Hand-written notes, 1972/73.

\bibitem{Juhasz}
I.~Juh{\'a}sz.
\newblock {\em Cardinal functions in topology}.
\newblock Mathematisch Centrum, Amsterdam, 1971.
\newblock In collaboration with A. Verbeek and N. S. Kroonenberg, Mathematical
  Centre Tracts, No. 34.

\bibitem{Koszmider}
Piotr Koszmider.
\newblock On the existence of strong chains in {$\frak{P}(\omega\sb 1)/{\rm
  Fin}$}.
\newblock {\em J. Symbolic Logic}, 63(3):1055--1062, 1998.

\bibitem{Koszmider2000}
Piotr Koszmider.
\newblock On strong chains of uncountable functions.
\newblock {\em Israel J. Math.}, 118:289--315, 2000.

\bibitem{Martinez}
Juan~Carlos Mart{\'{\i}}nez.
\newblock A consistency result on thin-very tall {B}oolean algebras.
\newblock {\em Israel J. Math.}, 123:273--284, 2001.

\bibitem{Morgan1996}
Charles Morgan.
\newblock Morasses, square and forcing axioms.
\newblock {\em Ann. Pure Appl. Logic}, 80(2):139--163, 1996.

\bibitem{Morgan1998}
Charles Morgan.
\newblock Higher gap morasses. {IA}. {G}ap-two morasses and condensation.
\newblock {\em J. Symbolic Logic}, 63(3):753--787, 1998.

\bibitem{Morgan2006}
Charles Morgan.
\newblock Local connectedness and distance functions.
\newblock In {\em Set theory}, Trends Math., pages 345--400. Birkh\"auser,
  Basel, 2006.

\bibitem{ShelahStanley}
Saharon Shelah and Lee Stanley.
\newblock A theorem and some consistency results in partition calculus.
\newblock {\em Ann. Pure Appl. Logic}, 36(2):119--152, 1987.

\bibitem{SolovayTennenbaum}
R.~M. Solovay and S.~Tennenbaum.
\newblock Iterated {C}ohen extensions and {S}ouslin's problem.
\newblock {\em Ann. of Math. (2)}, 94:201--245, 1971.

\bibitem{Stanley2}
Lee Stanley.
\newblock L-like models of set theory: Forcing, combinatorial principles, and
  morasses.
\newblock Dissertation, UC Berkeley, 1977.

\bibitem{Stanley}
Lee Stanley.
\newblock A short course on gap-one morasses with a review of the fine
  structure of {$L$}.
\newblock In {\em Surveys in set theory}, volume~87 of {\em London Math. Soc.
  Lecture Note Ser.}, pages 197--243. Cambridge Univ. Press, Cambridge, 1983.

\bibitem{Tennenbaum}
S.~Tennenbaum.
\newblock Souslin's problem.
\newblock {\em Proc. Nat. Acad. Sci. U.S.A.}, 59:60--63, 1968.

\bibitem{Todorcevic1985}
Stevo Todor{\v{c}}evi{\'c}.
\newblock Directed sets and cofinal types.
\newblock {\em Trans. Amer. Math. Soc.}, 290(2):711--723, 1985.

\bibitem{Todorcevic}
Stevo Todor{\v{c}}evi{\'c}.
\newblock Partitioning pairs of countable ordinals.
\newblock {\em Acta Math.}, 159(3-4):261--294, 1987.

\bibitem{Todorcevic1989}
Stevo Todor{\v{c}}evi{\'c}.
\newblock {\em Partition problems in topology}, volume~84 of {\em Contemporary
  Mathematics}.
\newblock American Mathematical Society, Providence, RI, 1989.

\bibitem{StevoBook}
Stevo Todorcevic.
\newblock {\em Walks on ordinals and their characteristics}, volume 263 of {\em
  Progress in Mathematics}.
\newblock Birkh\"auser Verlag, Basel, 2007.

\bibitem{Velleman1984}
Dan Velleman.
\newblock Simplified morasses.
\newblock {\em J. Symbolic Logic}, 49(1):257--271, 1984.

\bibitem{Velleman1987b}
Dan Velleman.
\newblock Gap-{$2$} morasses of height {$\omega$}.
\newblock {\em J. Symbolic Logic}, 52(4):928--938, 1987.

\bibitem{Velleman1987a}
Dan Velleman.
\newblock Simplified gap-{$2$} morasses.
\newblock {\em Ann. Pure Appl. Logic}, 34(2):171--208, 1987.

\end{thebibliography}

\end{document}